\documentclass[12pt]{article}
\usepackage{amsmath,amsfonts,amssymb}
\usepackage{theorem}
\theorembodyfont{\rmfamily}
\newtheorem{theorem}{Theorem}[section]
\newtheorem{proposition}[theorem]{Proposition}
\newtheorem{lemma}[theorem]{Lemma}
\newtheorem{definition}[theorem]{Definition}

\newtheorem{example}[theorem]{Example}
\newtheorem{remark}[theorem]{Remark}

\begin{document}
$\,$\vspace{10mm}

\begin{center}
{\textsf {\huge Crystal Interpretation of Kerov--Kirillov--}}
\vspace{2mm}\\
{\textsf {\huge Reshetikhin Bijection II}}
\vspace{6mm}\\
{\textsf {\LARGE --- Proof for $\mathfrak{sl}_{\mbox{n}}$ Case ---}}
\vspace{10mm}\\
{\textsf {\Large Reiho Sakamoto}}
\vspace{2mm}\\
{\textsf {Department of Physics, Graduate School of Science,}}
\vspace{-1mm}\\
{\textsf {University of Tokyo, Hongo, Bunkyo-ku,}}
\vspace{-1mm}\\
{\textsf {Tokyo, 113-0033, Japan}}
\vspace{-1mm}\\
{\textsf {reiho@monet.phys.s.u-tokyo.ac.jp}}
\vspace{15mm}
\end{center}

\begin{abstract}
\noindent
In proving the Fermionic formulae, combinatorial bijection called
the Kerov--Kirillov--Reshetikhin (KKR) bijection plays the central role.
It is a bijection between the set of highest paths and
the set of rigged configurations.
In this paper, we give a proof of crystal theoretic reformulation
of the KKR bijection.
It is the main claim of Part I written by
A. Kuniba, M. Okado, T. Takagi, Y. Yamada, and the author.
The proof is given by introducing a structure of
affine combinatorial $R$ matrices
on rigged configurations.
\end{abstract}
\pagebreak

\section{Introduction}
In this paper, we treat the relationship between the Fermionic
formulae and the well-known soliton cellular automata
``box-ball system.''
The Fermionic formulae are certain combinatorial identities,
and a typical example can be found
in the context of solvable lattice models.
The basis of these formulae is a combinatorial bijection called
the Kerov--Kirillov--Reshetikhin (KKR) bijection \cite{KKR,KR,KSS},
which gives one-to-one correspondences between the two
combinatorial objects called rigged configurations and
highest paths.
Precise description of the bijection is given in Section \ref{sec:defKKR}.

{}From the physical point of view, rigged configurations give an index set
for eigenvectors and eigenvalues of the Hamiltonian
that appears when we use the Bethe ansatz under the string hypothesis
(see, e.g., \cite{Sch1} for an introductory account of it),
and highest paths give an index set that appears when we use the
corner transfer matrix method (see, e.g., \cite{Bax}).
Therefore the KKR bijection means that
although neither the Bethe ansatz nor the corner transfer matrix method
is a rigorous mathematical theory, two index sets have
one-to-one correspondence.

Eventually, it becomes clearer that the KKR bijection itself possesses
a rich structure, especially with respect to the
representation theory of crystal bases \cite{Kas}.
For example, an extension of the rigged configuration called
unrestricted rigged configuration is recently introduced \cite{Sch2,DS},
and its crystal structure, i.e.,
actions of the Kashiwara operators on them is explicitly
determined.
It gives a natural generalization of the KKR bijection
which covers nonhighest weight elements.
(See, e.g., \cite{HKOTY1,HKOTT,SS}
for other information).

On the other hand, the box-ball system has entirely different background.
This model is a typical example of
soliton cellular automata introduced by Takahashi--Satsuma \cite{TS,Tak}.
It is an integrable discrete dynamical system
and has a direct connection with the discrete analogue
of the Lotka--Volterra equation
\cite{TTMS} (see also \cite{TTS}).
Though the time evolution of the system is described by a
simple combinatorial procedure, it beautifully exhibits a soliton dynamics.
Recently, a remarkable correspondence between the box-ball
systems and the crystal bases theory was discovered,
and it caused a lot of interests
(see, e.g., \cite{HIK,HKT,FOY,HHIKTT,HKOTY2} for related topics).
\vspace{4mm}

In Part I \cite{KOSTY} of this pair of papers, a unified treatment of both
the Fermionic formula (or the KKR bijection)
and the box-ball systems was presented.
It can be viewed as the inverse scattering
formalism (or Gelfand--Levitan formalism) for the box-ball systems.
In Part I, generalizations to
arbitrary nonexceptional affine Lie algebras
(the Okado--Schilling--Shimozono bijection \cite{OSS})
are also discussed.

In this paper, we give a proof of the result announced in Part I
for the general $\mathfrak{sl}_n$ case
(see Section 2.6 ``Main theorem'' of \cite{KOSTY}).
The precise statement of the result is formulated in
Theorem \ref{clcl-zdon} of Section \ref{sec:results} below.
According to our result,
the KKR bijection is interpreted in terms of combinatorial
$R$ matrices and energy functions of the crystals
(see Section \ref{sec:NYrule} for definitions).
Originally the KKR bijection is defined in a purely combinatorial way,
and it has no representation theoretic interpretation for a long time.
Therefore it is expected that our algebraic reformulation
will give some new insights into the theory
of crystals for finite-dimensional representations
of quantum affine Lie algebras \cite{KMN1,KMN2,KMN}.

Recently, as an application of our Theorem \ref{clcl-zdon},
explicit piecewise linear formula of the KKR bijection is
derived \cite{KSY}.
This formula involves the so-called tau functions which originate
from the theory of solitons \cite{JM}.
Interestingly, these tau functions have direct connection
with the Fermionic formula itself.
These results reveal unexpected link between the Fermionic formulae and
the soliton theory and, at the same time,
also give rise to general solution to the box-ball systems.
\vspace{4mm}

Let us describe some more details of our results.
As we have described before, main combinatorial objects
concerning the KKR bijection are rigged configurations
and highest paths.
Rigged configurations are the following set of data
\begin{equation}
\mathrm{RC}=\left( (\mu_i^{(0)}),\,
(\mu_i^{(1)},r_i^{(1)}),\,
\cdots ,
(\mu_i^{(n-1)},r_i^{(n-1)})\right) ,
\label{intro:level0}
\end{equation}
where $\mu^{(a)}_i\in\mathbb{Z}_{>0}$ and
$r^{(a)}_i\in\mathbb{Z}_{\geq 0}$ for $0\leq a\leq n-1$
and $1\leq i\leq l^{(a)}$ ($l^{(a)}\in\mathbb{Z}_{\geq 0}$).
They obey certain selection rule, which will be given in
Definition \ref{def:rc}.
On the other hand, highest paths are the highest weight elements
of $B_{k_1}\otimes B_{k_2}\otimes\cdots\otimes B_{k_N}$, where
$B_{k_i}$ is the crystal of $k_i$th symmetric power of the vector
(or natural) representation of
the quantum enveloping algebra $U_q(\mathfrak{sl}_n)$.
We regard elements of $B_{k_i}$ as row-type semi-standard
Young tableaux filled in with $k_i$ letters from
1 to $n$.
In this paper, we only treat a map from rigged configurations
to highest paths.

In order to reformulate the KKR bijection algebraically,
we notice that the nested structure arising on
rigged configuration Eq.(\ref{intro:level0}) is important.
More precisely, we introduce the following family
of subsets of RC for $0\leq a\leq n-1$;
\begin{equation}
\mathrm{RC}^{(a)}=
\left( (\mu_i^{(a)}),\,
(\mu_i^{(a+1)},r_i^{(a+1)}),\,
\cdots ,
(\mu_i^{(n-1)},r_i^{(n-1)})\right) .
\label{intro:levela}
\end{equation}
On this RC${}^{(a)}$, we can also apply the KKR bijection.
Then we obtain a path whose tensor factors are represented
by tableaux filled in with letters from 1 to $n-a$.
However, for our construction, it is convenient to
add $a$ to each letter contained in the path.
Thus, we assume that the path obtained from
RC${}^{(a)}$ contains letters $a+1$ to $n$.
Let us tentatively denote the resulting path $p^{(a)}$.
Then we can define the following maps:
\begin{equation}
p^{(a)}\xrightarrow{\Phi^{(a)}\circ C^{(a)}}
p^{(a-1)}.
\end{equation}
We postpone a presice definition of these maps
$\Phi^{(a)}\circ C^{(a)}$ until Section \ref{sec:mainresult},
but it should be stressed that the definition uses
only conbinatorial $R$ matrices and energy functions.
Note that the KKR bijection on RC${}^{(n-1)}$
trivially yields a path of the form
$p^{(n-1)}=\bigotimes_{i=1}^{l^{(n-1)}}
\fbox{$n^{\mu^{(n-1)}_i}$}$,
where \fbox{$n^\mu$} is a tableaux representation of crystals.
Therefore, by successive applications of
$\Phi^{(a)}C^{(a)}$ onto $p^{(n-1)}$,
we obtain the construction
\begin{equation}
p=
\Phi^{(1)}{C}^{(1)}
\Phi^{(2)}{C}^{(2)}\cdots
\Phi^{(n-1)}{C}^{(n-1)}
\left(
\bigotimes_{i=1}^{l^{(n-1)}}
\fbox{$n^{\mu^{(n-1)}_i}$}\right) ,
\end{equation}
where $p$ is the path corresponding to the original
rigged configuration RC (Eq.(\ref{intro:level0})).
\vspace{4mm}

The plan of this paper is as follows.
In Section \ref{sec:preliminaries}, we review definitions
of rigged configurations and the KKR bijection.
In Section \ref{sec:results}, we review combinatorial $R$
matrices and energy function following the graphical
rule in terms of winding and unwinding pairs
introduced in \cite{NY}.
We then define scattering data in Eqs.(\ref{def:sd}) and (\ref{defmode})
and define the operators $C^{(a)}$ and $\Phi^{(a)}$.
Our main result is formulated in Theorem \ref{clcl-zdon}.
The rest of the paper is devoted to a proof of this theorem.
In Section \ref{sec:normal}, we recall the
Kirillov--Schilling--Shimozono's result (Theorem \ref{KSS}).
This theorem describes the dependence of a resulting path
with respect to orderings of $\mu^{(0)}$ of RC.
We then introduce an important modification of rigged configurations.
More precisely, we replace $\mu^{(a)}$ of RC${}^{(a)}$
by $\mu^{(a)}\cup\mu^{(a+1)}\cup (1^L)$,
where the integer $L$ will be determined by Proposition \ref{Q-Q}.
We then apply Theorem \ref{KSS} to this modified rigged
configuration and obtain the isomorphism of Proposition \ref{moritsuke1}.
This reduces our remaining task to giving interpretation of modes
$d_i$ (Eq.(\ref{defmode})) in terms of the KKR bijection.
Example of these arguments is given in Example \ref{ex1}.
In Section \ref{sec:mode}, we connect modes $d_i$
with rigged configuration in Proposition \ref{Q-Q}.
By using this proposition, we introduce a structure
related with the energy function in Section \ref{sec:energy&KKR}.
This is described in Theorem \ref{unwinding}
(see also Examples \ref{ex:unwinding1} and
\ref{ex:unwinding2} as to the meanings of this theorem).
In Section \ref{sec:proof}, we give a proof of Theorem \ref{unwinding}
and hence complete a proof of Theorem \ref{clcl-zdon}.
We do this by directly connecting the graphical rule
of energy function given in Section \ref{sec:NYrule}
with rigged configuration.
In fact, we explicitly construct a structure of unwinding pairs on
the rigged configurations in Proposition \ref{NYpair}.

\section{Preliminaries}\label{sec:preliminaries}
\subsection{Rigged configurations}
In this section, we briefly review the Kerov--Kirillov--Reshetikhin
(KKR) bijection.
The KKR bijection gives one-to-one correspondences between the
set of rigged configurations and the set of highest weight elements
in tensor products of crystals of symmetric powers of the
vector (or natural) representation of $U_q(\mathfrak{sl}_n)$,
which we call paths.

Let us define the rigged configurations.
Consider the following collection of data:
\begin{equation}
\mu^{(a)}=
\left( \mu^{(a)}_1,\mu^{(a)}_2,\cdots ,\mu^{(a)}_{l^{(a)}}\right),
\qquad (0\leq a\leq n-1,\,
l^{(a)}\in\mathbb{Z}_{\geq 0},\,
\mu^{(a)}_i\in \mathbb{Z}_{>0}).
\end{equation}
We use usual Young diagrammatic expression for these integer sequences
$\mu^{(a)}$, although our $\mu^{(a)}$ are not necessarily
monotonically decreasing sequences.

\begin{definition}
$(1)$ For a given diagram $\mu$, we introduce coordinates
(row, column) of each boxes just like matrix entries.
For a box $\alpha$ of $\mu$, $col(\alpha)$ is column coordinate of $\alpha$.
Then we define the following subsets:
\begin{eqnarray}
\mu |_{\leq j}&:=&\{ \alpha|\alpha\in\mu\, ,col(\alpha )\leq j\},\\
\mu |_{>j}&:=&\{ \alpha|\alpha\in\mu\, ,col(\alpha )>j\}\, .
\end{eqnarray}

$(2)$ For a sequence of diagrams
$(\mu^{(0)},\mu^{(1)},\cdots ,\mu^{(n-1)})$,
we define $Q_j^{(a)}$ by 
\begin{equation}
Q_j^{(a)}:=\sum_{k=1}^{l^{(a)}}\min (j,\mu^{(a)}_k),
\end{equation}
i.e., the number of boxes in $\mu^{(a)}|_{\leq j}$.
Then the vacancy number $p^{(a)}_j$ for rows of
$\mu^{(a)}$ is defined by 
\begin{equation}
p^{(a)}_j:=Q^{(a-1)}_j-2Q^{(a)}_j+Q^{(a+1)}_j,
\end{equation}
where $j$ is the width of the corresponding row.
\rule{5pt}{10pt}
\end{definition}

\begin{definition}\label{def:rc}
Consider the following set of data:
\begin{equation}
\mathrm{RC}:=\left( (\mu_i^{(0)}),\,
(\mu_i^{(1)},r_i^{(1)}),\,
\cdots ,
(\mu_i^{(n-1)},r_i^{(n-1)})\right) .
\end{equation}

$(1)$ If all vacancy numbers for
$(\mu^{(1)},\mu^{(2)},\cdots ,\mu^{(n-1)})$
are nonnegative,
\begin{equation}
0\leq p^{(a)}_{\mu^{(a)}_i},\qquad
(1\leq a\leq n-1, 1\leq i\leq l^{(a)}),
\end{equation}
then RC is called a configuration.

$(2)$ If an integer $r^{(a)}_i$ satisfies the condition
\begin{equation}
0\leq r^{(a)}_i\leq p^{(a)}_{\mu^{(a)}_i},
\end{equation}
then $r^{(a)}_i$ is called a rigging associated with row
$\mu^{(a)}_i$.
For the rows of equal widths, i.e., $\mu^{(a)}_i=\mu^{(a)}_{i+1}$,
we assume that $r^{(a)}_i\leq r^{(a)}_{i+1}$.

$(3)$ If RC is a configuration and if all integers $r^{(a)}_i$
are riggings associated with row $\mu^{(a)}_i$,
then RC is called $\mathfrak{sl}_n$ rigged configuration.
\rule{5pt}{10pt}
\end{definition}
In the rigged configuration, $\mu^{(0)}$ is sometimes called a
quantum space which determines the shape of the corresponding path,
as we will see in the next subsection.
In the definition of the KKR bijection,
the following notion is important.
\begin{definition}
For a given rigged configuration,
consider a row $\mu^{(a)}_i$ and corresponding rigging $r^{(a)}_i$.
If they satisfy the condition
\begin{equation}
r^{(a)}_i=p^{(a)}_{\mu^{(a)}_i},
\end{equation}
then the row $\mu^{(a)}_i$ is called singular.
\rule{5pt}{10pt}
\end{definition}

\subsection{The KKR bijection}\label{sec:defKKR}
In this subsection, we define the KKR bijection.
In what follows, we treat a bijection $\phi$ to obtain a highest path $p$
from a given rigged configuration RC,
\begin{equation}
\phi :{\rm RC} \longrightarrow p\in
B_{k_N}\otimes\cdots\otimes B_{k_2}\otimes B_{k_1},
\end{equation}
where
\begin{equation}
{\rm RC}=
\left( (\mu_i^{(0)}),\,
(\mu_i^{(1)},r_i^{(1)}),\,
\cdots ,
(\mu_i^{(n-1)},r_i^{(n-1)})\right) ,
\end{equation}
is the  rigged configuration
defined in the last subsection, and
$N(=l^{(0)})$ is the length of the partition $\mu^{(0)}$.
$B_k$ is the crystal of the $k$th symmetric power of the vector (or natural)
representation of $U_q(\mathfrak{sl}_n)$.
As a set, it is equal to
\begin{equation}
B_k=\{(x_1,x_2,\cdots ,x_n)\in\mathbb{Z}^n_{\geq 0}\,\vert\,
x_1+x_2+\cdots +x_n=k\}.
\end{equation}
We usually identify elements of $B_k$ as the 
semi-standard Young tableaux
\begin{equation}
(x_1,x_2,\cdots ,x_n)=
\fbox{$\overbrace{1\cdots 1}^{x_1}\overbrace{2\cdots 2}^{x_2}
\cdots\cdots\overbrace{n\mathstrut\cdots n}^{x_n}$}\, ,
\end{equation}
i.e., the number of letters $i$ contained in a tableau is $x_i$.
\vspace{4mm}

\begin{definition}\label{defKKR}
For a given RC, the image (or path) $p$ of the KKR bijection $\phi$
is obtained by the following procedure.
\vspace{4mm}

\noindent
{\bf Step 1:}
For each row of the quantum space $\mu^{(0)}$, we re-assign the indices
from 1 to $N$ arbitrarily
and reorder it as the composition
\begin{equation}
\mu^{(0)}=\left( \mu^{(0)}_N,\cdots ,\mu^{(0)}_2, \mu^{(0)}_1\right) .
\end{equation}
Take the row $\mu^{(0)}_1$.
Recall that $\mu^{(0)}$ is not necessarily
monotonically decreasing integer sequence.

\noindent
{\bf Step 2:}
We denote each box of the row $\mu^{(0)}_1$ as follows:
\begin{equation}
\mu^{(0)}_1=
\unitlength 23pt
\begin{picture}(5,1)(-0.3,0.3)
\multiput(0,0)(0,1){2}{\line(1,0){5}}
\multiput(0,0)(1,0){2}{\line(0,1){1}}
\put(0.1,0.3){$\alpha^{(0)}_{l_1}$}
\multiput(1.35,0.35)(0.3,0){5}{$\cdot$}
\multiput(3,0)(1,0){3}{\line(0,1){1}}
\put(3.1,0.3){$\alpha^{(0)}_{2}$}
\put(4.1,0.3){$\alpha^{(0)}_{1}$}
\put(5.1,0.3){.}
\end{picture}
\end{equation}
Corresponding to the row $\mu^{(0)}_1$, we take $p_1$
as the following array of $l_1$ empty boxes:
\begin{equation}
p_1\, =
\unitlength 23pt
\begin{picture}(5,1)(-0.3,0.3)
\multiput(0,0)(0,1){2}{\line(1,0){5}}
\multiput(0,0)(1,0){2}{\line(0,1){1}}
\multiput(1.35,0.35)(0.3,0){5}{$\cdot$}
\multiput(3,0)(1,0){3}{\line(0,1){1}}
\put(5.1,0.3){.}
\end{picture}
\end{equation}
Starting from the box $\alpha^{(0)}_1$, we recursively take
$\alpha^{(i)}_1\in\mu^{(i)}$ by the following Rule 1.
\begin{quotation}
\noindent
{\bf Rule 1:}
Assume that we have already chosen $\alpha^{(i-1)}_1\in\mu^{(i-1)}$.
Let $g^{(i)}$ be the set of all rows of $\mu^{(i)}$ whose
widths $w$ satisfy
\begin{equation}
w\geq col(\alpha^{(i-1)}_1).
\end{equation}
Let $g^{(i)}_s$ ($\subset g^{(i)}$) be the set of all singular
rows (i.e., its rigging is equal to the vacancy
number of the corresponding row)
in a set $g^{(i)}$.
If $g^{(i)}_s\neq\emptyset$, then choose one of the shortest rows
of $g^{(i)}_s$ and denote by $\alpha^{(i)}_1$ its rightmost box.
If $g^{(i)}_s=\emptyset$, then we take
$\alpha^{(i)}_1=$ $\cdots$ $=\alpha^{(n-1)}_1$ $=\emptyset$.
\end{quotation}

\noindent
{\bf Step 3:}
{}From RC
remove the boxes $\alpha^{(0)}_1$, $\alpha^{(1)}_1$, $\cdots$,
$\alpha^{(j_1-1)}_1$ chosen above,
where $j_1-1$ is defined by
\begin{equation}
j_1-1=\max_{0\leq k\leq n-1,\atop \alpha^{(k)}_1\neq\emptyset} k.
\end{equation}
After removal, the new RC is obtained by the following Rule 2.
\begin{quotation}
\noindent
{\bf Rule 2:}
Calculate again all the vacancy numbers
$p^{(a)}_i=Q^{(a-1)}_i-2Q^{(a)}_i+Q^{(a+1)}_i$ according to the removed RC.
For a row which is not removed,
take the rigging equal to the corresponding
rigging before removal.
For a row which is removed, take the rigging equal to the new
vacancy number of the corresponding row.
\end{quotation}
Put the letter $j_1$ into the leftmost empty box of $p_1$:
\begin{equation}
p_1\, =
\unitlength 23pt
\begin{picture}(5,1)(-0.3,0.3)
\multiput(0,0)(0,1){2}{\line(1,0){5}}
\multiput(0,0)(1,0){3}{\line(0,1){1}}
\put(0.3,0.4){$j_1$}
\multiput(2.35,0.35)(0.3,0){5}{$\cdot$}
\multiput(4,0)(1,0){2}{\line(0,1){1}}
\put(5.1,0.3){.}
\end{picture}
\end{equation}

\noindent
{\bf Step 4:}
Repeat Step 2 and Step 3 for the rest of boxes $\alpha^{(0)}_2$,
$\alpha^{(0)}_3$, $\cdots$, $\alpha^{(0)}_{l_1}$
in this order.
Put the letters $j_2,j_3,\cdots ,j_{l_1}$ into empty boxes of $p_1$ from left to right.

\noindent
{\bf Step 5:}
Repeat Step 1 to Step 4 for the rest of rows $\mu^{(0)}_2$,
$\mu^{(0)}_3$, $\cdots$, $\mu^{(0)}_N$ in this order.
Then we obtain $p_k$ from $\mu^{(0)}_k$,
which we identify with the element of $B_{\mu^{(0)}_k}$.
Then we obtain
\begin{equation}
p=p_N\otimes\cdots\otimes p_2\otimes p_1
\end{equation}
as an image of $\phi$.
\rule{5pt}{10pt}
\end{definition}
Note that the resulting image $p$ is a function of the ordering of
$\mu^{(0)}$ which we choose in Step 1.
Its dependence is described in Theorem \ref{KSS} below.

The above procedure is summarized in the following diagram.
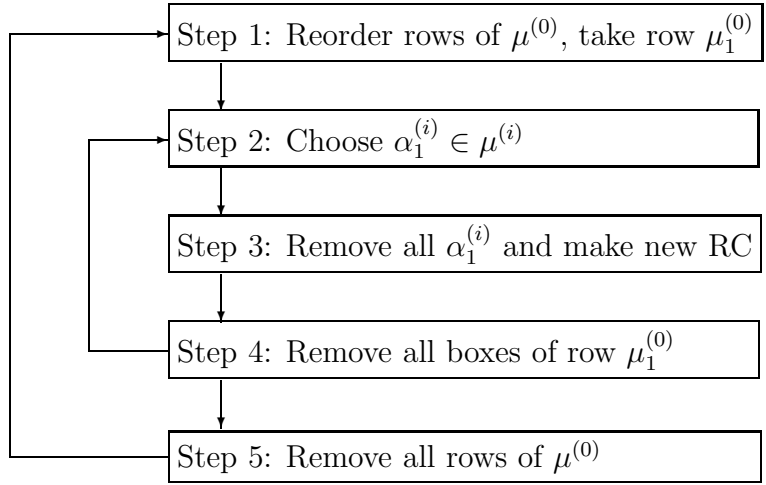
\begin{center}
\unitlength 20pt
\begin{picture}(15,9)
\put(3,8){\fbox{Step 1: Reorder rows of $\mu^{(0)}$, take row $\mu^{(0)}_1$}}
\put(3,6){\fbox{Step 2: Choose $\alpha^{(i)}_1\in\mu^{(i)}$\hspace{30mm}}}
\put(3,4){\fbox{Step 3: Remove all $\alpha^{(i)}_1$ and make new RC}}
\put(3,2){\fbox{Step 4: Remove all boxes of row $\mu^{(0)}_1$\hspace{10mm}}}
\put(3,0){\fbox{Step 5: Remove all rows of $\mu^{(0)}$\hspace{20mm}}}
\multiput(4,1.65)(0,2){4}{\vector(0,-1){0.9}}
\put(3,0.2){\line(-1,0){3}}
\put(0,0.2){\line(0,1){8}}
\put(0,8.2){\vector(1,0){3}}
\put(3,2.2){\line(-1,0){1.5}}
\put(1.5,2.2){\line(0,1){4}}
\put(1.5,6.2){\vector(1,0){1.5}}
\end{picture}
\end{center}

\begin{example}\label{ex:kkr}
We give one simple but nontrivial example.
Consider the following $\mathfrak{sl}_3$ rigged configuration:
\begin{center}
\unitlength 13pt
\begin{picture}(11,5.5)
\multiput(0,0)(1,0){2}{\line(0,1){4}}
\multiput(0,0)(0,1){3}{\line(1,0){1}}
\multiput(0,3)(0,1){2}{\line(1,0){2}}
\put(2,3){\line(0,1){1}}
\put(0.5,4.5){$\mu^{(0)}$}
\put(5,2){\line(1,0){1}}
\multiput(5,2)(1,0){2}{\line(0,1){2}}
\multiput(5,3)(0,1){2}{\line(1,0){2}}
\put(7,3){\line(0,1){1}}
\put(5.5,4.5){$\mu^{(1)}$}
\put(4.2,2.15){1}
\put(4.2,3.15){0}
\put(6.3,2.15){0}
\put(7.3,3.15){0}
\multiput(10,3)(1,0){2}{\line(0,1){1}}
\multiput(10,3)(0,1){2}{\line(1,0){1}}
\put(9.2,3.15){0}
\put(11.3,3.15){0}
\put(10,4.5){$\mu^{(2)}$}
\end{picture}
\end{center}
We write the vacancy number on the left and riggings on the right
of the Young diagrams.
We reorder $\mu^{(0)}$ as $( 1,1,2,1)$;
thus we remove the following boxes
\unitlength 13pt
\begin{picture}(1,1)(0,0.2)
\multiput(0,0)(1,0){2}{\line(0,1){1}}
\multiput(0,0)(0,1){2}{\line(1,0){1}}
\put(0.15,0.25){$\times$}
\end{picture} :
\begin{center}
\unitlength 13pt
\begin{picture}(11,5.5)
\multiput(0,0)(1,0){2}{\line(0,1){4}}
\multiput(0,0)(0,1){3}{\line(1,0){1}}
\multiput(0,3)(0,1){2}{\line(1,0){2}}
\put(2,3){\line(0,1){1}}
\put(0.5,4.5){$\mu^{(0)}$}
\put(0.15,2.25){$\times$}
\put(5,2){\line(1,0){1}}
\multiput(5,2)(1,0){2}{\line(0,1){2}}
\multiput(5,3)(0,1){2}{\line(1,0){2}}
\put(7,3){\line(0,1){1}}
\put(5.5,4.5){$\mu^{(1)}$}
\put(4.2,2.15){1}
\put(4.2,3.15){0}
\put(6.3,2.15){0}
\put(7.3,3.15){0}
\put(6.15,3.25){$\times$}
\multiput(10,3)(1,0){2}{\line(0,1){1}}
\multiput(10,3)(0,1){2}{\line(1,0){1}}
\put(9.2,3.15){0}
\put(11.3,3.15){0}
\put(10,4.5){$\mu^{(2)}$}
\end{picture}
\end{center}
We obtain $p_1=\fbox{2}$.
Note that, in this step, we cannot remove singular
row of $\mu^{(2)}$, since it is shorter than 2.

After removing two boxes, calculate again the vacancy numbers
and make the row of $\mu^{(1)}$ (which is removed) singular.
Then we obtain the following configuration:
\begin{center}
\unitlength 13pt
\begin{picture}(11,4.5)
\multiput(0,0)(1,0){2}{\line(0,1){3}}
\multiput(0,0)(0,1){2}{\line(1,0){1}}
\multiput(0,2)(0,1){2}{\line(1,0){2}}
\put(2,2){\line(0,1){1}}
\put(0.5,3.5){$\mu^{(0)}$}
\put(1.15,2.25){$\times$}
\multiput(5,1)(1,0){2}{\line(0,1){2}}
\multiput(5,1)(0,1){3}{\line(1,0){1}}
\multiput(4.2,1.15)(0,1){2}{0}
\multiput(6.3,1.15)(0,1){2}{0}
\put(5,3.5){$\mu^{(1)}$}
\multiput(9,2)(1,0){2}{\line(0,1){1}}
\multiput(9,2)(0,1){2}{\line(1,0){1}}
\put(8.2,2.15){0}
\put(10.3,2.15){0}
\put(9,3.5){$\mu^{(2)}$}
\end{picture}
\end{center}
Next, we remove the box
\unitlength 13pt
\begin{picture}(1,1)(0,0.2)
\multiput(0,0)(1,0){2}{\line(0,1){1}}
\multiput(0,0)(0,1){2}{\line(1,0){1}}
\put(0.15,0.25){$\times$}
\end{picture}
from the above configuration.
We cannot remove $\mu^{(1)}$, since all singular rows
are shorter than 2.
Thus, we obtain
$p_2=\unitlength 13pt
\begin{picture}(2,1)(0,0.2)
\multiput(0,0)(1,0){3}{\line(0,1){1}}
\multiput(0,0)(0,1){2}{\line(1,0){2}}
\put(0.3,0.2){1}
\end{picture}$
, and the new rigged configuration is the following:
\begin{center}
\unitlength 13pt
\begin{picture}(10,4.5)
\multiput(0,0)(1,0){2}{\line(0,1){3}}
\multiput(0,0)(0,1){4}{\line(1,0){1}}
\put(0,3.5){$\mu^{(0)}$}
\put(0.15,2.25){$\times$}
\multiput(4,1)(1,0){2}{\line(0,1){2}}
\multiput(4,1)(0,1){3}{\line(1,0){1}}
\multiput(3.2,1.15)(0,1){2}{0}
\multiput(5.3,1.15)(0,1){2}{0}
\put(4,3.5){$\mu^{(1)}$}
\put(4.15,2.25){$\times$}
\multiput(8,2)(1,0){2}{\line(0,1){1}}
\multiput(8,2)(0,1){2}{\line(1,0){1}}
\put(7.2,2.15){0}
\put(9.3,2.15){0}
\put(8,3.5){$\mu^{(2)}$}
\put(8.15,2.25){$\times$}
\end{picture}
\end{center}

This time, we can remove $\mu^{(1)}$ and $\mu^{(2)}$
and obtain $p_2=\unitlength 13pt
\begin{picture}(2,1)(0,0.2)
\multiput(0,0)(1,0){3}{\line(0,1){1}}
\multiput(0,0)(0,1){2}{\line(1,0){2}}
\put(0.3,0.2){1}
\put(1.3,0.2){3}
\end{picture}$ .
Then we obtain the following configuration:
\begin{center}
\unitlength 13pt
\begin{picture}(9,3.5)
\multiput(0,0)(1,0){2}{\line(0,1){2}}
\multiput(0,0)(0,1){3}{\line(1,0){1}}
\put(0,2.5){$\mu^{(0)}$}
\put(0.15,1.25){$\times$}
\multiput(4,1)(1,0){2}{\line(0,1){1}}
\multiput(4,1)(0,1){2}{\line(1,0){1}}
\put(4,2.5){$\mu^{(1)}$}
\multiput(3.2,1.15)(0,1){1}{0}
\multiput(5.3,1.15)(0,1){1}{0}
\put(4.15,1.25){$\times$}
\put(8,1.25){$\emptyset$}
\put(8,2.5){$\mu^{(2)}$}
\end{picture}
\end{center}
{}From this configuration we remove the boxes
\unitlength 13pt
\begin{picture}(1,1)(0,0.2)
\multiput(0,0)(1,0){2}{\line(0,1){1}}
\multiput(0,0)(0,1){2}{\line(1,0){1}}
\put(0.15,0.25){$\times$}
\end{picture}
and obtain $p_3=\fbox{2}$,
and the new configuration becomes the following:
\begin{center}
\unitlength 13pt
\begin{picture}(9,2.5)
\multiput(0,0)(1,0){2}{\line(0,1){1}}
\multiput(0,0)(0,1){2}{\line(1,0){1}}
\put(0,1.5){$\mu^{(0)}$}
\put(0.15,0.25){$\times$}
\put(4,0.25){$\emptyset$}
\put(4,1.5){$\mu^{(1)}$}
\put(8,0.25){$\emptyset$}
\put(8,1.5){$\mu^{(2)}$}
\end{picture}
\end{center}
Finally we obtain $p_4=\fbox{1}$.

To summarize, we obtain
\begin{equation}
p=\fbox{1}\otimes\fbox{2}\otimes\fbox{13}
\otimes\fbox{2}\, ,
\end{equation}
as an image of the KKR bijection.
\rule{5pt}{10pt}

\end{example}

\section{Crystal base theory and the KKR bijection}\label{sec:results}
\subsection{Combinatorial $R$ matrix and energy functions}\label{sec:NYrule}
In this section, we formulate the statement of our main result.
First of all, let us summarize the basic objects from the crystal bases
theory, namely, the combinatorial $R$ matrix and associated energy function.

For two crystals $B_k$ and $B_l$ of $U_q(\mathfrak{sl}_n)$,
one can define the tensor product
$B_k\otimes B_l=\{b\otimes b'\mid b\in B_k,b'\in B_l\}$.
Then we have a unique isomorphism $R:B_k\otimes B_l
\stackrel{\sim}{\rightarrow}B_l\otimes B_k$, i.e. a unique map
which commutes with actions of the Kashiwara operators.
We call this map combinatorial $R$ matrix
and usually write the map $R$ simply by $\simeq$.

Following Rule 3.11 of \cite{NY}, we introduce a graphical rule
to calculate the combinatorial $R$ matrix for $\mathfrak{sl}_n$ 
and the energy function.
Given the two elements
\[
x=(x_1,x_2,\cdots ,x_{n})\in B_k, \quad 
y=(y_1,y_2,\cdots ,y_{n})\in B_l,
\]
we draw the following diagram to represent the tensor
product $x\otimes y$:
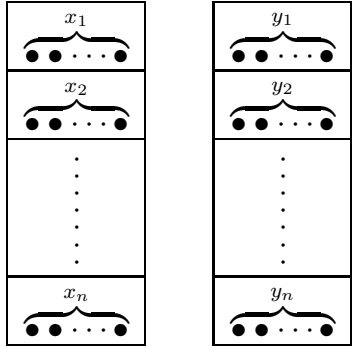
\begin{center}
\unitlength 13pt
\begin{picture}(10,10.5)
\multiput(0,0)(6,0){2}{
\multiput(0,0)(4,0){2}{\line(0,1){10}}
\multiput(0,0)(0,2){2}{\line(1,0){4}}
\multiput(0,6)(0,2){3}{\line(1,0){4}}
}
\put(0.5,0.2){$\overbrace{\bullet\bullet\cdots\bullet}^{x_{n}}$}
\put(0.5,6.2){$\overbrace{\bullet\bullet\cdots\bullet}^{x_2}$}
\put(0.5,8.2){$\overbrace{\bullet\bullet\cdots\bullet}^{x_1}$}
\multiput(1.9,2.2)(0,0,5){7}{$\cdot$}
\put(6.5,0.2){$\overbrace{\bullet\bullet\cdots\bullet}^{y_{n}}$}
\put(6.5,6.2){$\overbrace{\bullet\bullet\cdots\bullet}^{y_2}$}
\put(6.5,8.2){$\overbrace{\bullet\bullet\cdots\bullet}^{y_1}$}
\multiput(7.9,2.2)(0,0,5){7}{$\cdot$}
\end{picture}
\end{center}
The combinatorial $R$ matrix and energy function $H$ for
$B_k\otimes B_l$ (with $k\geq l$) are calculated by
the following rule.
\begin{enumerate}
\item
Pick any dot, say $\bullet_a$, in the right column and connect it
with a dot $\bullet_a'$ in the left column by a line.
The partner $\bullet_a'$ is chosen
{}from the dots which are in the lowest row among all dots
whose positions are higher than that of $\bullet_a$.
If there is no such a dot, we return to the bottom, and
the partner $\bullet_a'$ is chosen from the dots
in the lowest row among all dots.
In the former case, we call such a pair ``unwinding,"
and, in the latter case, we call it ``winding."

\item
Repeat procedure (1) for the remaining unconnected dots
$(l-1)$ times.

\item
Action of the combinatorial $R$ matrix is obtained by
moving all unpaired dots in the left column to the right
horizontally.
We do not touch the paired dots during this move.

\item
The energy function $H$ is given by the number of winding pairs.
\end{enumerate}

The number of winding (or unwinding) pairs is sometimes called
the winding (or unwinding, respectively) number of tensor product.
It is known that the resulting combinatorial $R$ matrix
and the energy functions are not affected by the
order of making pairs
(\cite{NY}, Propositions 3.15 and 3.17).
For more properties, including that the above
definition indeed satisfies the axiom, see \cite{NY}.

\begin{example}
The diagram for $\fbox{1344}\otimes\fbox{234}$ is
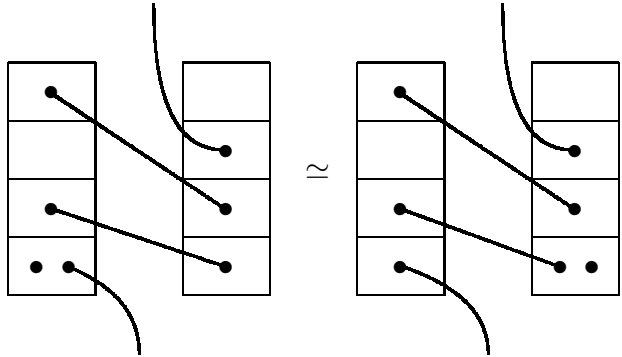
\begin{center}
\unitlength 11pt
\begin{picture}(10,14)(6,-1)
\multiput(0,2)(6,0){2}{
\multiput(0,0)(3,0){2}{\line(0,1){8}}
\multiput(0,0)(0,2){5}{\line(1,0){3}}}
\put(1.2,8.7){$\bullet$}
\put(1.2,4.7){$\bullet$}
\put(0.7,2.7){$\bullet$}
\put(1.8,2.7){$\bullet$}
\put(7.2,6.7){$\bullet$}
\put(7.2,4.7){$\bullet$}
\put(7.2,2.7){$\bullet$}
\thicklines
\qbezier(5,12)(5,7)(7.4,7)
\qbezier(1.4,9.0)(1.4,9.0)(7.4,5)
\qbezier(1.4,5.0)(1.4,5.0)(7.4,3)
\qbezier(2,3)(4.5,2)(4.5,0)
\thinlines
\put(10.2,6){$\simeq$}
\put(12,0){
\multiput(0,2)(6,0){2}{
\multiput(0,0)(3,0){2}{\line(0,1){8}}
\multiput(0,0)(0,2){5}{\line(1,0){3}}}
\put(1.2,8.7){$\bullet$}
\put(1.2,4.7){$\bullet$}
\put(1.2,2.7){$\bullet$}
\put(7.2,6.7){$\bullet$}
\put(7.2,4.7){$\bullet$}
\put(6.7,2.7){$\bullet$}
\put(7.8,2.7){$\bullet$}
\thicklines
\qbezier(5,12)(5,7)(7.4,7)
\qbezier(1.4,9.0)(1.4,9.0)(7.4,5)
\qbezier(1.4,5.0)(1.4,5.0)(6.9,3.0)
\qbezier(1.5,3)(4.5,2)(4.5,0)
}
\end{picture}
\end{center}

By moving the unpaired dot (letter 4) in the left column to 
the right, we obtain
\begin{equation*}
\fbox{1344}\otimes\fbox{234}
\simeq 
\fbox{134}\otimes\fbox{2344}\, .
\end{equation*}
Since we have one winding pair and two unwinding pairs,
the energy function is
$H\left(
\fbox{1344}\otimes\fbox{234}
\right)=1$.
\rule{5pt}{10pt}
\end{example}
By the definition, the winding numbers for 
$x \otimes y$ and ${\tilde y} \otimes {\tilde x}$ are the same 
if $x \otimes y \simeq {\tilde y} \otimes {\tilde x}$ by the 
combinatorial $R$ matrix.

\subsection{Formulation of the main result}\label{sec:mainresult}
{}From now on, we reformulate the original KKR bijection
in terms of the combinatorial $R$ and energy function.
Consider the $\mathfrak{sl}_n$ rigged configuration
as follows:
\begin{equation}
{\rm RC}=
\left( (\mu_i^{(0)}),\,
(\mu_i^{(1)},r_i^{(1)}),\,
\cdots ,
(\mu_i^{(n-1)},r_i^{(n-1)})\right) .
\label{level0}
\end{equation}
By applying the KKR bijection, we obtain a path $\tilde{s}^{(0)}$.

In order to obtain a path $\tilde{s}^{(0)}$ by algebraic procedure,
we have to introduce a nested structure on the rigged configuration.
More precisely, we consider the following subsets
of given configuration (\ref{level0}) for $0\leq a\leq n-1$:
\begin{equation}
\mathrm{RC}^{(a)}:=
\left( (\mu_i^{(a)}),\,
(\mu_i^{(a+1)},r_i^{(a+1)}),\,
\cdots ,
(\mu_i^{(n-1)},r_i^{(n-1)})\right) .
\label{levela}
\end{equation}
RC${}^{(a)}$ is a $\mathfrak{sl}_{n-a}$ rigged configuration,
and RC${}^{(0)}$ is nothing but the original RC.
Therefore we can perform the KKR bijection on RC${}^{(a)}$
and obtain a path $\tilde{s}^{(a)}$ with letters $1,2,\cdots ,n-a$.
However, for our construction, it is convenient to add $a$ to
all letters in a path.
Thus we assume that a path $\tilde{s}^{(a)}$ contains
letters $a+1,\cdots ,n$.

As in the original path $\tilde{s}^{(0)}$,
we should consider $\tilde{s}^{(a)}$ as highest weight elements
of tensor products of crystals as follows:
\begin{equation}
\tilde{s}^{(a)}=b_1\otimes\cdots\otimes b_N\in
B_{k_1}\otimes\cdots\otimes
B_{k_N}, \; \; (k_i = \mu^{(a)}_i,\, N=l^{(a)}).
\end{equation}
The meaning of crystals $B_k$ here is as follows.
$B_k$ is crystal of the $k$th symmetric power representation
of the vector (or natural) representation of $U_q(\mathfrak{sl}_{n-a})$.
As a set, it is equal to
\begin{equation}
B_k=\{(x_{a+1},x_{a+2},\cdots ,x_n)\in\mathbb{Z}^{n-a}_{\geq 0}\,\vert\,
x_{a+1}+x_{a+2}+\cdots +x_n=k\}.
\end{equation}
We can identify elements of $B_k$ as 
semi-standard Young tableaux containing letters $a+1,\cdots ,n$.
Also, we can naturally extend the graphical rule for the combinatorial $R$
matrix and energy function (see Section \ref{sec:NYrule}) to this case.
The highest weight element of $B_k$ takes the form
\begin{equation}
\fbox{$(a+1)^k$}=\fbox{$(a+1)\cdots (a+1)^{\mathstrut}$}\in B_k.
\end{equation}
This corresponds to the so-called lower diagonal embedding of
$\mathfrak{sl}_{n-a}$ into $\mathfrak{sl}_{n}$.

{}From now on,
let us construct an element of affine crystal $s^{(a)}$
from $\tilde{s}^{(a)}$ combined with information of
riggings $r^{(a)}_i$,
\begin{equation}\label{def:sd}
s^{(a)}:=
b_1[d_1]\otimes
\cdots\otimes
b_N[d_N]\, \in
{\rm aff}(B_{k_1})\otimes
\cdots\otimes{\rm aff}(B_{k_N}).
\end{equation}
Here $\mathrm{aff}(B)$ is the affinization of a crystal $B$.
As a set, it is equal to
\begin{equation}
\mathrm{aff}(B)=\{b[d]\, |\, d\in\mathbb{Z},b\in B\},
\end{equation}
where integers $d$ of $b[d]$ are often called modes.
We can extend the combinatorial $R$:
$B\otimes B'\simeq B'\otimes B$ to the affine case 
$\mathrm{aff}(B)\otimes\mathrm{aff}(B')
\simeq\mathrm{aff}(B')\otimes\mathrm{aff}(B)$
by the relation
\begin{equation}
b[d]\otimes b'[d']\simeq
\tilde{b}'[d'-H(b\otimes b')]\otimes \tilde{b}[d+H(b\otimes b')],
\label{eq:R}
\end{equation}
where $b\otimes b'\simeq\tilde{b}'\otimes \tilde{b}$
is the isomorphism of combinatorial $R$ matrix for classical crystals
which was defined in Section \ref{sec:NYrule}.

Now we define the element $s^{(a)}$ of Eq.(\ref{def:sd})
from a path $\tilde{s}^{(a)}$ and riggings $r^{(a)}_i$.
Mode $d_i$ of $b_i[d_i]$ of $s^{(a)}$ is defined by the formula
\begin{equation}\label{defmode}
d_i:=r^{(a)}_i+\sum_{0\leq l< i}
H\left( b_l\otimes b^{(l+1)}_i
\right) ,\quad b_0:= \fbox{$(a+1)^{\max k_i}$}\, ,
\end{equation}
where $r^{(a)}_i$ is the rigging corresponding to a row $\mu^{(a)}_i$
of RC${}^{(0)}$ which yielded the element $b_i$ of $\tilde{s}^{(a)}$.
The elements $b^{(l+1)}_i$ $(l<i)$ are defined by sending $b_i$ successively
to the right of $b_l$ under the
isomorphism of combinatorial $R$ matrices:
\begin{eqnarray}
&&
b_1\otimes\cdots\otimes
b_l\otimes b_{l+1}\otimes\cdots\otimes
b_{i-2}\otimes b_{i-1}\otimes
b_i\otimes\cdots\nonumber\\
&\simeq&
b_1\otimes\cdots\otimes
b_l\otimes b_{l+1}\otimes\cdots\otimes
b_{i-2}\otimes b_{i}^{(i-1)}\otimes
b_{i-1}'\otimes\cdots\nonumber\\
&\simeq&\cdots\cdots\nonumber\\
&\simeq&
b_1\otimes
\cdots\otimes
b_l\otimes
b^{(l+1)}_i\otimes\cdots\otimes
b_{i-3}'\otimes b_{i-2}'\otimes
b_{i-1}'\otimes\cdots .
\end{eqnarray}
This definition of $d_i$ is compatible with the following
commutation relation of affine combinatorial $R$ matrix:
\begin{equation}
\cdots\otimes b_i[d_i]\otimes b_{i+1}[d_{i+1}]\otimes\cdots
\,\simeq\,
\cdots\otimes b_{i+1}'[d_{i+1}-H]\otimes
b_i'[d_i+H]\otimes\cdots,
\end{equation}
where $b_i\otimes b_{i+1}\simeq b_{i+1}'\otimes b_i'$ is an isomorphism
by classical combinatorial $R$ matrix (see Theorem \ref{KSS} below)
and $H=H(b_i\otimes b_{i+1})$.
We call an element of affine crystal $s^{(a)}$ a scattering data.

For a scattering data $s^{(a)}=b_1[d_1]\otimes\cdots\otimes b_N[d_N]$
obtained from the quantum space $\mu^{(a)}$,
we define the normal ordering as follows.
\begin{definition}\label{norder}
For a given scattering data $s^{(a)}$, we define
the sequence of subsets
\begin{equation}
\mathcal{S}_1\subset
\mathcal{S}_2\subset\cdots\subset
\mathcal{S}_N\subset
\mathcal{S}_{N+1}\,
\end{equation}
as follows.
$\mathcal{S}_{N+1}$ is the set of all permutations
which are obtained by $\widehat{\mathfrak{sl}}_{n-a}$
combinatorial $R$ matrices acting
on each tensor product in $s^{(a)}$.
$\mathcal{S}_i$ is the subset of $\mathcal{S}_{i+1}$ consisting of
all the elements of $\mathcal{S}_{i+1}$
whose $i$th mode from the left end
are maximal in $\mathcal{S}_{i+1}$.
Then the elements of $\mathcal{S}_1$ are called the normal ordered form
of $s^{(a)}$.
\rule{5pt}{10pt}
\end{definition}

Although the above normal ordering is not unique,
we choose any one of the normal ordered scattering data
which is obtained from the path $\tilde{s}^{(a)}$ and denote it by
$C^{(a)}(\tilde{s}^{(a)})$.
See Remark \ref{characterization} for alternative
characterization of the normal ordering.
For $C^{(a)}(\tilde{s}^{(a)})=b_1[d_1]\otimes\cdots\otimes b_N[d_N]$
($b_i\in B_{k_i}$),
we define the following element of $\mathfrak{sl}_{n-a+1}$
crystal with letters $a,\cdots ,n$:
\begin{equation}\label{def:c}
c=
\fbox{$a$}^{\,\otimes d_1}\otimes
b_1\otimes
\fbox{$a$}^{\,\otimes (d_2-d_1)}\otimes
b_2\otimes\cdots\otimes
\fbox{$a$}^{\,\otimes (d_N-d_{N-1})}\otimes
b_N\, .
\end{equation}
In the following, we need the map $C^{(n-1)}$.
To define it, we use combinatorial $R$ of
``$\widehat{\mathfrak{sl}}_1$'' crystal defined as follows:
\begin{equation}
\fbox{$n^k$}_{\, d_2}\otimes\fbox{$n^l$}_{\, d_1}\simeq
\fbox{$n^l$}_{\, d_1-H}\otimes\fbox{$n^k$}_{\, d_2+H},
\end{equation}
where $H$ is now $H=\min (k,l)$, and we have denoted
$b_k[d_k]$ as $\fbox{$b_k$}_{\, d_k}$.
This is a special case of the combinatorial $R$ matrix and energy
function defined in Section \ref{sec:NYrule},
and $\widehat{\mathfrak{sl}}_1$ corresponds to
the $\mathfrak{sl}_2$ subalgebra generated by $e_0$ and $f_0$.

We introduce another operator $\Phi^{(a)}$,
\begin{equation}
\Phi^{(a)}:
{\rm aff}(B_{k_1})\otimes\cdots\otimes
{\rm aff}(B_{k_N}) \rightarrow
B_{l_1}\otimes \cdots \otimes B_{l_{N'}}
\end{equation}
where we denote $l_i=\mu^{(a-1)}_i$ and $N'=l^{(a-1)}$.
$\Phi^{(a)}$ is defined by the following isomorphism of
$\mathfrak{sl}_{n-a+1}$ combinatorial $R$:
\begin{equation}\label{extra1}
\Phi^{(a)}\left( C^{(a)}(\tilde{s}^{(a)})\right)\otimes\left(
\bigotimes_{i=1}^{N}
\fbox{$a^{k_i}$}\right)
\otimes\fbox{$a$}^{\,\otimes d_N}
\,\simeq\,
c\otimes
\left(
\bigotimes_{i=1}^{N'}
\fbox{$a^{l_i}$}
\right) ,
\end{equation}
where $c$ is defined in Eq.(\ref{def:c}).

Then our main result is the following:
\begin{theorem}\label{clcl-zdon}
For the rigged configuration RC${}^{(a)}$
(see Eq.(\ref{levela})),
we consider the KKR bijection with
letters from $a+1$ to $n$ .
Then its image is given by
\begin{equation}
\Phi^{(a+1)}{C}^{(a+1)}
\Phi^{(a+2)}{C}^{(a+2)}\cdots
\Phi^{(n-1)}{C}^{(n-1)}
\left(
\bigotimes^{l^{(n-1)}}_{i=1}
\fbox{$n^{\mu^{(n-1)}_i}$}\right) .
\end{equation}
In particular, the KKR image $p$ of rigged configuration
(\ref{level0}) satisfies
\begin{equation}
p=
\Phi^{(1)}{C}^{(1)}
\Phi^{(2)}{C}^{(2)}\cdots
\Phi^{(n-1)}{C}^{(n-1)}
\left(
\bigotimes^{l^{(n-1)}}_{i=1}
\fbox{$n^{\mu^{(n-1)}_i}$}\right) .
\end{equation}
The image of this map is independent of the choice
of maps $C^{(a)}$.
\rule{5pt}{10pt}
\end{theorem}

In practical calculation of this procedure,
it is convenient to introduce the following diagrams.
First, we express the isomorphism of the combinatorial
$R$ matrix
\begin{equation}
a\otimes b\simeq b^{'}\otimes a^{'}
\end{equation}
by the following vertex diagram:
\begin{center}
\unitlength 13pt
\begin{picture}(4,4)
\put(0.6,2.0){\line(1,0){2}}
\put(1.6,1.0){\line(0,1){2}}
\put(0,1.8){$a$}
\put(1.4,0){$b^{'}$}
\put(1.4,3.2){$b$}
\put(2.8,1.8){$a^{'}$}
\put(3.5,1.7){.}
\end{picture}
\end{center}
If we apply combinatorial $R$ successively as
\begin{equation}
a\otimes b\otimes c\simeq b^{'}\otimes a^{'}\otimes c
\simeq b^{'}\otimes c^{'}\otimes a^{''},
\end{equation}
then we express this by joining two vertices as follows:
\begin{center}
\unitlength 13pt
\begin{picture}(8,4)
\multiput(0,0)(3.2,0){2}{
\put(0.6,2.0){\line(1,0){2}}
\put(1.6,1.0){\line(0,1){2}}
}
\put(-0.1,1.8){$a$}
\put(1.4,0){$b^{'}$}
\put(1.4,3.2){$b$}
\put(2.9,1.8){$a^{'}$}
\put(4.6,3.2){$c$}
\put(4.6,0){$c^{'}$}
\put(6.1,1.8){$a^{''}$}
\put(7.0,1.7){.}
\end{picture}
\end{center}
Also, it is sometimes convenient to use the
notation $a\stackrel{H}{\otimes}b$
if we have $H=H(a\otimes b)$.

\begin{example}
We give an example of Theorem \ref{clcl-zdon}
along with the same rigged configuration
we have considered in Example \ref{ex:kkr}.
\begin{center}
\unitlength 13pt
\begin{picture}(11,5.5)
\multiput(0,0)(1,0){2}{\line(0,1){4}}
\multiput(0,0)(0,1){3}{\line(1,0){1}}
\multiput(0,3)(0,1){2}{\line(1,0){2}}
\put(2,3){\line(0,1){1}}
\put(0.5,4.5){$\mu^{(0)}$}
\put(5,2){\line(1,0){1}}
\multiput(5,2)(1,0){2}{\line(0,1){2}}
\multiput(5,3)(0,1){2}{\line(1,0){2}}
\put(7,3){\line(0,1){1}}
\put(5.5,4.5){$\mu^{(1)}$}
\put(4.2,2.15){1}
\put(4.2,3.15){0}
\put(6.3,2.15){0}
\put(7.3,3.15){0}
\multiput(10,3)(1,0){2}{\line(0,1){1}}
\multiput(10,3)(0,1){2}{\line(1,0){1}}
\put(9.2,3.15){0}
\put(11.3,3.15){0}
\put(10,4.5){$\mu^{(2)}$}
\end{picture}
\end{center}

First we calculate a path $\tilde{s}^{(2)}$,
which is an image of the following rigged configuration
(it contains the quantum space only):
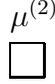
\begin{center}
\unitlength 13pt
\begin{picture}(11,2.5)(5,3)
\multiput(10,3)(1,0){2}{\line(0,1){1}}
\multiput(10,3)(0,1){2}{\line(1,0){1}}
\put(10,4.5){$\mu^{(2)}$}
\end{picture}
\end{center}
The KKR bijection trivially yields its image as
\begin{equation}
\tilde{s}^{(2)}=\fbox{3}\, .
\nonumber
\end{equation}
We define the mode of \fbox{3} using Eq.(\ref{defmode}).
We put $b_0=\fbox{3}$ and $b_1=\fbox{3}\, (=\tilde{s}^{(2)})$.
Since we have $\fbox{3}\stackrel{1}{\otimes}\fbox{3}$
and $r^{(2)}_1=0$, the mode is $0+1=1$.
Therefore we have
\begin{equation}
C^{(2)}\left(\fbox{3}\right)
=\fbox{3}_{\, 1}.
\nonumber
\end{equation}
Note that $\fbox{3}_{\, 1}$ is trivially normal ordered.

Next we calculate $\Phi^{(2)}$.
Let us take the numbering of rows of $\mu^{(1)}$
as $(\mu^{(1)}_1,\mu^{(1)}_2)=(2,1)$,
i.e., the resulting path is an element of
$B_{\mu^{(1)}_1}\otimes B_{\mu^{(1)}_2}=B_2\otimes B_1$.
{}From $\fbox{3}_{\, 1}$, we create an element $\fbox{2}\otimes\fbox{3}$
(see Eq.(\ref{def:c}))
and consider the following tensor product
(see the right-hand side of Eq.(\ref{extra1})):
\begin{equation}\nonumber
\fbox{2}\otimes\fbox{3}\otimes
\left(\fbox{22}\otimes\fbox{2}
\right) .
\end{equation}
We move \fbox{3} to the right of
$\fbox{22}\otimes\fbox{2}$ and next we
move \fbox{2} to the right, as in the following diagram:
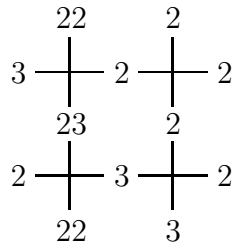
\begin{center}
\unitlength 13pt
\begin{equation}\nonumber
\begin{picture}(6,7)
\multiput(0,1)(0,3){2}{
\multiput(0,0)(3,0){2}{
\put(0,1){\line(1,0){2}}
\put(1,0){\line(0,1){2}}
}}
\put(-0.7,1.7){2}
\put(-0.7,4.7){3}
\put(0.6,0.1){22}
\put(0.6,3.2){23}
\put(0.6,6.3){22}
\put(2.3,1.7){3}
\put(2.3,4.7){2}
\put(3.8,0.1){3}
\put(3.8,3.2){2}
\put(3.8,6.3){2}
\put(5.3,1.7){2}
\put(5.3,4.7){2}
\end{picture}
\end{equation}
\end{center}
We have omitted framings of tableaux \fbox{$\ast$}
in the above diagram.
Therefore we have
\begin{equation}
\Phi^{(2)}\left(
\fbox{3}_{\, 1}
\right)
=\fbox{22}\otimes\fbox{3}\, .
\nonumber
\end{equation}
Note that the result depend on the choice of
the shape of path ($B_2\otimes B_1$).

Let us calculate $C^{(1)}$.
First, we determine the modes $d_1$, $d_2$ of
$\fbox{22}_{\, d_1}\otimes\fbox{3}_{\, d_2}$.
For $d_1$, we put $b_0=\fbox{22}$ , and the corresponding
value of an energy function is
$\fbox{22}\stackrel{2}{\otimes}\fbox{22}\otimes\fbox{3}$ ,
and the rigging is $r^{(1)}_1=0$;
hence we have $d_1=2+0=2$.
For $d_2$, we need the following values of
energy functions;
$\fbox{22}\otimes\fbox{22}\stackrel{0}{\otimes}\fbox{3}
\stackrel{R}{\simeq}
\fbox{22}\stackrel{1}{\otimes}\fbox{2}\otimes\fbox{23}$,
and the rigging is $r^{(1)}_2=0$.
Hence we have $d_2=0+1+0=1$.
In order to determine the normal ordering of
$\fbox{22}_{\, 2}\otimes\fbox{3}_{\, 1}$
($\stackrel{R}{\simeq}\fbox{2}_{\, 1}\otimes\fbox{23}_{\, 2}$),
following Definition \ref{norder},
we construct the set $\mathcal{S}_3$ as
\begin{equation}
\mathcal{S}_3=\left\{
\fbox{22}_{\, 2}\otimes\fbox{3}_{\, 1},
\fbox{2}_{\, 1}\otimes\fbox{23}_{\, 2}
\right\}.
\nonumber
\end{equation}
Therefore the normal ordered form is
\begin{equation}
C^{(1)}\left(
\fbox{22}\otimes\fbox{3}
\right) =
\fbox{2}_{\, 1}\otimes\fbox{23}_{\, 2}.
\nonumber
\end{equation}

Finally, we calculate $\Phi^{(1)}$.
We assume that the resulting path is an element of
$B_1\otimes B_1\otimes B_2\otimes B_1$.
{}From $\fbox{2}_{\, 1}\otimes\fbox{23}_{\, 2}$
we construct an element
$\fbox{1}\otimes\fbox{2}\otimes\fbox{1}\otimes\fbox{23}$.
We consider the tensor product
\begin{equation}\label{ex:clclzdon_eq}
\fbox{1}\otimes\fbox{2}\otimes\fbox{1}\otimes\fbox{23}\otimes
\left(
\fbox{1}\otimes \fbox{1}\otimes \fbox{11}\otimes \fbox{1}
\right)
\end{equation}
and apply combinatorial $R$ matrices successively as follows:
\begin{center}
\unitlength 13pt
\begin{equation}\label{ex:clclzdon}
\begin{picture}(13,13)(-1,0)
\multiput(0,1)(0,3){4}{
\multiput(0,0)(3,0){4}{
\put(0,1){\line(1,0){2}}
\put(1,0){\line(0,1){2}}
}}
\put(-0.7,1.7){1}
\put(-0.7,4.7){2}
\put(-0.7,7.7){1}
\put(-1.1,10.7){23}
\put(0.8,0.1){1}
\put(0.8,3.2){2}
\put(0.8,6.2){1}
\put(0.8,9.2){3}
\put(0.8,12.2){1}
\put(2.2,1.7){2}
\put(2.2,4.7){1}
\put(2.2,7.7){3}
\put(2.0,10.7){12}
\put(3.8,0.1){2}
\put(3.8,3.2){1}
\put(3.8,6.2){3}
\put(3.8,9.2){2}
\put(3.8,12.2){1}
\put(5.3,1.7){1}
\put(5.3,4.7){3}
\put(5.3,7.7){2}
\put(5.1,10.7){11}
\put(6.6,0.1){13}
\put(6.6,3.2){23}
\put(6.6,6.2){12}
\put(6.6,9.2){11}
\put(6.6,12.2){11}
\put(8.2,1.7){2}
\put(8.2,4.7){1}
\put(8.2,7.7){1}
\put(8.0,10.7){11}
\put(9.8,0.1){2}
\put(9.8,3.2){1}
\put(9.8,6.2){1}
\put(9.8,9.2){1}
\put(9.8,12.2){1}
\put(11.2,1.7){1}
\put(11.2,4.7){1}
\put(11.2,7.7){1}
\put(11.2,10.7){11}
\end{picture}
\end{equation}
\end{center}

Hence we obtain a path $\fbox{1}\otimes\fbox{2}\otimes
\fbox{13}\otimes\fbox{2}$,
which reconstructs a calculation of
Example \ref{ex:kkr}.
\rule{5pt}{10pt}
\end{example}
\begin{remark}
In the above calculation of $\Phi^{(2)}$,
we have assumed the shape of path as $B_2\otimes B_1$.
Then we calculated modes and obtained
$\fbox{22}_{\, 2}\otimes\fbox{3}_{\, 1}$.
Now suppose the path of the form $B_1\otimes B_2$ on the contrary.
In this case, calculation proceeds as follows:
\begin{center}
\unitlength 13pt
\begin{picture}(6,7)
\multiput(0,1)(0,3){2}{
\multiput(0,0)(3,0){2}{
\put(0,1){\line(1,0){2}}
\put(1,0){\line(0,1){2}}
}}
\put(-0.7,1.7){2}
\put(-0.7,4.7){3}
\put(0.8,0.1){2}
\put(0.8,3.2){3}
\put(0.8,6.3){2}
\put(2.3,1.7){3}
\put(2.3,4.7){2}
\put(3.6,0.1){23}
\put(3.6,3.2){22}
\put(3.6,6.3){22}
\put(5.3,1.7){2}
\put(5.3,4.7){2}
\end{picture}
\end{center}
{}From the values of energy functions
$\fbox{22}\stackrel{1}{\otimes}\fbox{2}\otimes\fbox{23}$
and
$\fbox{22}\otimes\fbox{2}\stackrel{0}{\otimes}\fbox{23}
\stackrel{R}{\simeq}
\fbox{22}\stackrel{2}{\otimes}\fbox{22}\otimes\fbox{3}$
and the riggings
$r^{(1)}_1=r^{(1)}_2=0$
we obtain an element
$\fbox{2}_{\, 1}\otimes\fbox{23}_{\, 2}$.
Comparing both results, we have
\begin{equation}
\fbox{2}_{\, 1}\otimes\fbox{23}_{\, 2}
\stackrel{R}{\simeq}
\fbox{22}_{\, 2}\otimes\fbox{3}_{\, 1}.
\nonumber
\end{equation}
This is a general consequence of the definition of mode (Eq.(\ref{defmode}))
and Theorem \ref{KSS} below.
\rule{5pt}{10pt}
\end{remark}

The rest of this paper is devoted to a proof of Theorem \ref{clcl-zdon}.

\section{Normal ordering from the KKR bijection}\label{sec:normal}
In the rest of this paper, we adopt the following numbering for
factors of the scattering data:
\begin{equation}
b_N[d_N]\otimes\cdots\otimes b_2[d_2]\otimes b_1[d_1]\,
\in{\rm aff}(B_{k_N})\otimes
\cdots\otimes{\rm aff}(B_{k_2})\otimes{\rm aff}(B_{k_1}),
\end{equation}
since this is more convenient when we are discussing about the relation between
the scattering data and KKR bijection.

It is known that the KKR bijection on rigged configuration RC
admits a structure of the combinatorial $R$ matrices.
This is described by the following powerful theorem
proved by Kirillov, Schilling and Shimozono (Lemma 8.5 of \cite{KSS}),
which plays an important role in the subsequent discussion.
\begin{theorem}\label{KSS}
Pick out any two rows from the quantum space $\mu ^{(0)}$ and
denote these by $\mu_a$ and $\mu_b$.
When we remove $\mu_a$ at first
and next $\mu_b$ by the KKR bijection,
then we obtain tableaux $\mu_a$ and $\mu_b$ with letters
$1,\cdots ,n$, which we denote by $A_1$ and $B_1$, respectively.
Next, on the contrary, we first remove $\mu_b$
and second $\mu_a$ (keeping the order of other removal invariant)
and we get $B_2$ and $A_2$.
Then we have
\begin{equation}
B_1\otimes A_1\,\simeq\, A_2\otimes B_2 ,
\end{equation}
under the isomorphism of $\mathfrak{sl}_n$ combinatorial $R$ matrix.
\rule{5pt}{10pt}
\end{theorem}

Our first task is to interpret the normal ordering which appear
in Definition \ref{norder}
in terms of purely KKR language.
We can achieve this translation if we make some tricky modification
on the rigged configuration.
Consider the rigged configuration
\begin{equation}
\mathrm{RC}^{(a-1)}=
\left( (\mu_i^{(a-1)}),\,
(\mu_i^{(a)},r_i^{(a)}),\,
\cdots ,
(\mu_i^{(n-1)},r_i^{(n-1)})\right) .
\label{levela-1}
\end{equation}
Then modify its quantum space $\mu^{(a-1)}$ as
\begin{equation}
\mu^{(a-1)}_+:=\mu^{(a-1)}\cup
\mu^{(a)}\cup (1^L)\, ,
\end{equation}
where $L$ is some sufficiently large integer to be determined below.
For the time being, we take $L$ large enough so that
configuration $\mu^{(a)}$
never becomes singular while we are removing $\mu^{(a-1)}$ part
from quantum space $\mu^{(a-1)}\cup \mu^{(a)}\cup (1^L)$
under the KKR procedure.
Then we obtain the modified rigged configuration
\begin{equation}
\mathrm{RC}^{(a-1)}_+:=
\left( (\mu_{+i}^{(a-1)}),\,
(\mu_i^{(a)},r_i^{(a)}),\,
\cdots ,
(\mu_i^{(n-1)},r_i^{(n-1)})\right) ,
\label{levela+}
\end{equation}
where $\mu_{+i}^{(a-1)}$ is the $i$th row of the quantum space
$\mu^{(a-1)}_+$.
In subsequent discussions, we always assume this modified
form of the quantum space unless otherwise stated.

For the KKR bijection on rigged configuration RC${}^{(a-1)}_+$,
we have two different ways to remove rows of quantum space
$\mu_+^{(a-1)}$.
We describe these two cases respectively.\\
{\bf Case 1.}
Remove $\mu^{(a)}$ and $(1^L)$ from $\mu_+^{(a-1)}$.
Then the rigged configuration RC${}^{(a-1)}_+$ reduces
to the original rigged configuration RC${}^{(a-1)}$.
Let us write the KKR image of RC${}^{(a-1)}$ by $p$,
then the KKR image of modified rigged configuration
RC${}^{(a-1)}_+$ in this case is
\begin{equation}
p\otimes\left(
\bigotimes_{i=1}^{l^{(a)}}\fbox{$a^{\mu^{(a)}_i}$}
\right)\otimes
\fbox{$a$}^{\,\otimes L}.
\end{equation}\\
{\bf Case 2.}
Remove $\mu^{(a-1)}$ from $\mu_+^{(a-1)}$ in RC${}^{(a-1)}_+$,
then quantum space becomes $\mu^{(a)}\cup (1^L)$.
Next, we remove the boxes of $(1^L)$ one by one until
some rows in $\mu^{(a)}$ become singular.
At this time, we choose any one of the singular rows in
$\mu^{(a)}$ and call it $\mu^{(a)}_1$.
We also define an integer $d_1(\leq L)$ such that
$(1^L)$ part of the quantum space is now
reduced to $(1^{d_1})$.
Then we have the following:
\begin{lemma}\label{lifting}
In the above setting, we remove row $\mu^{(a)}_1$ in
$\mu_+^{(a-1)}$ by the KKR procedure (with letters $a,\cdots ,n$)
and obtain a tableau $b_1\in B_{k_1}$.
On the other hand, consider the KKR procedure (with letters $a+1,\cdots ,n$)
on rigged configuration RC${}^{(a)}$, and remove row
$\mu^{(a)}_1$ as a first step of the procedure.
Then we obtain the same tableau $b_1$. 
\end{lemma}
{\bf Proof.}
Consider the rigged configuration RC${}^{(a-1)}_+$
after removing $\mu^{(a-1)}\cup(1^{L-d_1})$ from $\mu_+^{(a-1)}$.
When we begin to remove row $\mu^{(a)}_1$ in the quantum space
$\mu_+^{(a-1)}$,
we first remove the rightmost box of the row $\mu^{(a)}_1$, call box $x$.
Then, by the definition of $d_1$, the row $\mu^{(a)}_1$ in the next configuration
$\mu^{(a)}$ is singular so that we can remove the rightmost
box of the row $\mu^{(a)}_1\subset\mu^{(a)}$.
After removing $x$, the remaining row
$\mu^{(a)}_1\setminus\{ x\}\subset\mu^{(a)}$
is made to be singular again.

In the next step, we remove the box $x'\in\mu_+^{(a-1)}$ which is
on the left of the box $x$.
Then we can remove the corresponding box $x'\in\mu^{(a)}$.
Continuing in this way, we
remove both rows $\mu^{(a)}_1$ in quantum space
$\mu_+^{(a-1)}$ and $\mu^{(a)}$ simultaneously.
We see that this box removing operations on $\mu^{(a)}$,
$\mu^{(a+1)}$, $\cdots$, $\mu^{(n-1)}$ of RC${}^{(a-1)}_+$
coincides with the one that we have when we remove
$\mu^{(a)}_1$ of the quantum space of RC${}^{(a)}$.
\rule{5pt}{10pt}\vspace{4mm}

Let us return to the description of Case 2 procedure,
where we have just removed both $\mu^{(a)}_1$ from quantum space
$\mu^{(a-1)}_+$ and $\mu^{(a)}$.
Again, we remove boxes of $(1^{d_1})$ part of the quantum space
one by one until some singular rows appear in partition $\mu^{(a)}$
and choose any one of the singular rows, which we call $\mu^{(a)}_2$.
At this moment, the part $(1^{d_1})$ is reduced to $(1^{d_2})$.
We then remove both $\mu^{(a)}_2$ in quantum space and $\mu^{(a)}$
just as in the above lemma and obtain a tableau $b_2$.

We do this procedure recursively until all boxes of the quantum
space are removed.
Therefore the KKR image in this Case 2 is
\begin{equation}
\fbox{$a$}^{\,\otimes d_N}\otimes
b_N\otimes\cdots\otimes
\fbox{$a$}^{\,\otimes (d_2-d_3)}\otimes
b_2\otimes
\fbox{$a$}^{\,\otimes (d_1-d_2)}\otimes
b_1\otimes
\left(
\bigotimes^{l^{(a-1)}}_{i=1}\fbox{$a^{\mu^{(a-1)}_i}$}
\right) ,
\end{equation}
where we write $N=l^{(a)}$ and substitute $L$ in
$\mu^{(a-1)}_+$ by $d_1$.
This completes a description of Case 2 procedure.\vspace{4mm}

Note that, in this expression, the letter $a$ is separated from
the letters $a+1,\cdots ,n$ contained in $b_i$.
By virtue of this property, we introduce the following:
\begin{definition}\label{KKRnorder}
In the above Case 2 procedure, we obtain $b_i$ and the associated
integers $d_i$ by the KKR bijection.
{}From this data we construct the element
\begin{equation}
\tilde{C}^{(a)}
:=b_N[d_N]\otimes\cdots\otimes b_1[d_1]\in
{\rm aff}(B_{k_N})\otimes\cdots\otimes{\rm aff}(B_{k_1})
\end{equation}
and call this a KKR normal ordered product.
\rule{5pt}{10pt}
\end{definition}
To obtain a KKR normal ordering, we have to refer the actual KKR procedure.
Although the KKR normal ordering $\tilde{C}^{(a)}$ has not been identified
with the one defined in Definition \ref{norder},
these two procedures provide the interpretation of $\Phi^{(a)}$
operator.
More precisely, for each product
\begin{equation}
c:=\fbox{$a$}^{\,\otimes d_N}\otimes
b_N\otimes\cdots\otimes
\fbox{$a$}^{\,\otimes (d_2-d_3)}\otimes
b_2\otimes
\fbox{$a$}^{\,\otimes (d_1-d_2)}\otimes
b_1
\end{equation}
constructed from $\tilde{C}^{(a)}$,
we have the following isomorphism.
\begin{proposition}\label{moritsuke1}
For the rigged configuration RC${}^{(a-1)}_+$, we have
\begin{equation}
p\otimes\left(
\bigotimes_{i=1}^{l^{(a)}}\fbox{$a^{\mu^{(a)}_i}$}
\right)\otimes
\fbox{$a$}^{\,\otimes d_1}
\,\simeq\,
c\otimes\left(
\bigotimes^{l^{(a-1)}}_{i=1}\fbox{$a^{\mu^{(a-1)}_i}$}
\right) ,
\end{equation}
where the isomorphism is given by the $\mathfrak{sl}_{n-a+1}$ combinatorial $R$
matrix with letters $a,\cdots ,n$,
and $p$ is a path obtained by the KKR bijection on the original rigged
configuration RC${}^{(a-1)}$.
\end{proposition}
{\bf Proof.}
{}From the above construction we see that a difference between Case 1
and Case 2 is just the difference of order of removing rows of
$\mu^{(a-1)}_+$ in RC${}^{(a-1)}_+$.
Hence we can apply Theorem \ref{KSS} to claim that both
expressions are mutually isomorphic.
\rule{5pt}{10pt}\vspace{4mm}

This is just the $\Phi^{(a)}$ part of Theorem \ref{clcl-zdon}.
We continue to study further properties of this KKR normal ordered
product $\tilde{C}^{(a)}$.
Let us perform the above Case 2 procedure on RC${}^{(a-1)}_+$ and
obtain the KKR normal ordered product
\begin{equation}
b_{i_N}[d_{i_N}]\otimes\cdots\otimes
b_{i_2}[d_{i_2}]\otimes b_{i_1}[d_{i_1}],
\end{equation}
where each tableau $b_{i_k}$ comes from a row $\mu^{(a)}_{i_k}$.
However, there is an ambiguity in the choice of singular rows in Case 2.
Assume that we obtain another KKR normal ordered product
\begin{equation}
b_{j_N}'[d_{j_N}']\otimes\cdots\otimes
b_{j_2}'[d_{j_2}']\otimes b_{j_1}'[d_{j_1}']
\end{equation}
from the same configuration RC${}^{(a-1)}_+$.
We assume that each tableau $b_{j_k}'$ comes from a row $\mu^{(a)}_{j_k}$.
Then these two products have the following property.
\begin{lemma}
In this setting, we have
\begin{equation}\label{bbbb=bbbb}
b_{i_N}\otimes\cdots\otimes b_{i_2}\otimes b_{i_1}
\,\simeq\,
b_{j_N}'\otimes\cdots\otimes b_{j_2}'\otimes b_{j_1}'
\end{equation}
by $\mathfrak{sl}_{n-a}$ combinatorial $R$ matrices.
\end{lemma}
{\bf Proof.}
By using the argument of the proof of Lemma \ref{lifting},
we regard the left-hand side of Eq.(\ref{bbbb=bbbb})
as the path obtained from RC${}^{(a)}$ by removing
rows $\mu^{(a)}_{i_1},\cdots ,\mu^{(a)}_{i_N}$ of
the quantum space in this order.
Similarly, the right-hand side of Eq.(\ref{bbbb=bbbb})
is the path obtained by removing rows
$\mu^{(a)}_{j_1},\cdots ,\mu^{(a)}_{j_N}$ in this order.
Therefore we can apply Theorem \ref{KSS} to
obtain the isomorphism.
\rule{5pt}{10pt}\vspace{4mm}

\begin{example}\label{ex1}
We give an example of general argument given in this section
along with the following rigged configuration RC:
\begin{equation}
\begin{array}{ll}
(\mu^{(0)})&=\left( 1^{13}\right) ,\\
(\mu^{(1)},r^{(1)})&=\bigl( (4,0),(3,1),(1,4)\bigl),\\
(\mu^{(2)},r^{(2)})&=\bigl( (2,0),(1,0)\bigl),\\
(\mu^{(3)},r^{(3)})&=\bigl( (1,0)\bigl);
\end{array}\nonumber
\end{equation}
in the diagrammatic expression, it is
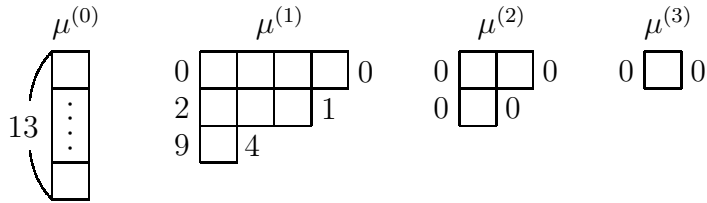
\begin{center}
\unitlength 14pt
\begin{picture}(20,5)
\multiput(1,0)(1,0){2}{\line(0,1){4}}
\multiput(1,0)(0,1){2}{\line(1,0){1}}
\multiput(1,3)(0,1){2}{\line(1,0){1}}
\multiput(1.37,1.18)(0,0.4){4}{$\cdot$}
\put(-0.2,1.7){13}
\qbezier(1,0)(0.5,0.5)(0.4,1.3)
\qbezier(1,4)(0.5,3.5)(0.4,2.7)
\put(1,4.5){$\mu^{(0)}$}
\put(5,1){\line(1,0){1}}
\multiput(5,1)(1,0){2}{\line(0,1){3}}
\put(5,2){\line(1,0){3}}
\multiput(5,3)(0,1){2}{\line(1,0){4}}
\multiput(7,2)(1,0){2}{\line(0,1){2}}
\put(9,3){\line(0,1){1}}
\put(4.3,1.2){9}
\put(4.3,2.2){2}
\put(4.3,3.2){0}
\put(6.2,1.17){4}
\put(8.27,2.17){1}
\put(9.26,3.19){0}
\put(6.5,4.5){$\mu^{(1)}$}
\put(12,2){\line(1,0){1}}
\multiput(12,2)(1,0){2}{\line(0,1){2}}
\multiput(12,3)(0,1){2}{\line(1,0){2}}
\put(14,3){\line(0,1){1}}
\put(11.3,2.2){0}
\put(11.3,3.2){0}
\put(13.25,2.2){0}
\put(14.25,3.2){0}
\put(12.5,4.5){$\mu^{(2)}$}
\multiput(17,3)(1,0){2}{\line(0,1){1}}
\multiput(17,3)(0,1){2}{\line(1,0){1}}
\put(16.3,3.2){0}
\put(18.25,3.2){0}
\put(17,4.5){$\mu^{(3)}$}
\end{picture}
\end{center}
For each Young diagram, we assign the vacancy numbers (on the left)
and riggings (on the right) of the corresponding rows
(for example, the vacancy numbers of $\mu^{(1)}$ are 0,2,9,
and the corresponding riggings are 0,1,4, respectively).
By the usual KKR bijection, we obtain the following image (path) $p$:
\begin{equation}
p=\fbox{1}\otimes\fbox{1}\otimes\fbox{1}\otimes\fbox{1}\otimes
\fbox{2}\otimes\fbox{2}\otimes\fbox{3}\otimes\fbox{2}\otimes\fbox{1}\otimes
\fbox{4}\otimes\fbox{3}\otimes\fbox{2}\otimes\fbox{2} .
\nonumber
\end{equation}

In the next section, we will obtain a formula which determines
the mode $d_1$ (Proposition \ref{Q-Q}).
Using the formula, we calculate $d_1$ as follows:
\begin{eqnarray}\label{555}
d_1&=&\max\{
Q^{(1)}_1-Q^{(2)}_1+r^{(1)}_1,\,
Q^{(1)}_3-Q^{(2)}_3+r^{(1)}_3,\,
Q^{(1)}_4-Q^{(2)}_4+r^{(1)}_4\}
\nonumber\\
&=&\max\{
3-2+4,\, 7-3+1,\, 8-3+0\}
\nonumber\\
&=&\max\{5,5,5\}=5.
\end{eqnarray}
Thus, in the modified rigged configuration RC${}^{(0)}_+$
(Eq.(\ref{levela+})),
we have to take a quantum space as follows:
\begin{equation}
\mu^{(0)}_+=\mu^{(1)}\cup (1^{13})\cup (1^5)
=\{4,3,1,1^{18}\} .
\end{equation}
The modified rigged configuration RC${}^{(0)}_+$ in this case
takes the following shape:
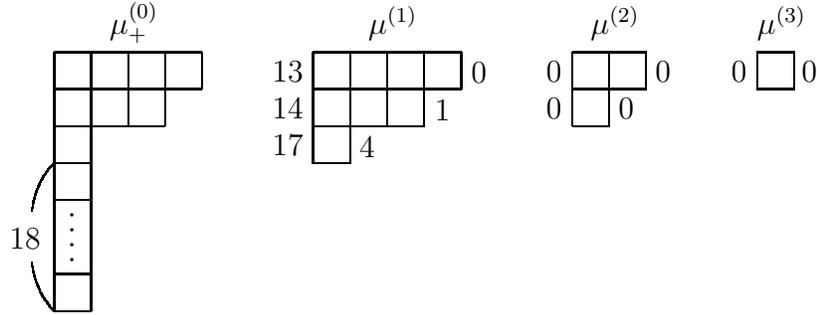
\begin{center}
\unitlength 14pt
\begin{picture}(22,8.5)
\multiput(1,0)(1,0){2}{\line(0,1){4}}
\multiput(1,0)(0,1){2}{\line(1,0){1}}
\multiput(1,3)(0,1){2}{\line(1,0){1}}
\multiput(1.37,1.18)(0,0.4){4}{$\cdot$}
\put(-0.2,1.7){18}
\qbezier(1,0)(0.5,0.5)(0.4,1.3)
\qbezier(1,4)(0.5,3.5)(0.4,2.7)
\multiput(1,4)(1,0){2}{\line(0,1){3}}
\put(1,5){\line(1,0){3}}
\multiput(1,6)(0,1){2}{\line(1,0){4}}
\multiput(3,5)(1,0){2}{\line(0,1){2}}
\put(5,6){\line(0,1){1}}
\put(2.5,7.5){$\mu^{(0)}_+$}
\put(8,4){\line(1,0){1}}
\multiput(8,4)(1,0){2}{\line(0,1){3}}
\put(8,5){\line(1,0){3}}
\multiput(8,6)(0,1){2}{\line(1,0){4}}
\multiput(10,5)(1,0){2}{\line(0,1){2}}
\put(12,6){\line(0,1){1}}
\put(6.9,4.2){17}
\put(6.9,5.2){14}
\put(6.9,6.2){13}
\put(9.25,4.17){4}
\put(11.29,5.17){1}
\put(12.28,6.17){0}
\put(9.5,7.5){$\mu^{(1)}$}
\put(15,5){\line(1,0){1}}
\multiput(15,5)(1,0){2}{\line(0,1){2}}
\multiput(15,6)(0,1){2}{\line(1,0){2}}
\put(17,6){\line(0,1){1}}
\put(14.3,5.2){0}
\put(14.3,6.2){0}
\put(16.25,5.2){0}
\put(17.25,6.2){0}
\put(15.5,7.5){$\mu^{(2)}$}
\multiput(20,6)(1,0){2}{\line(0,1){1}}
\multiput(20,6)(0,1){2}{\line(1,0){1}}
\put(19.3,6.2){0}
\put(21.2,6.2){0}
\put(20,7.5){$\mu^{(3)}$}
\end{picture}
\end{center}

We remove boxes according to Case 1 procedure given above.
In this procedure, we remove the $\mu^{(1)}\cup (1^5)$ part
from the quantum space $\mu^{(0)}_+$.
Then the remaining configuration is exactly equal to the original one
whose quantum space is $\mu^{(0)}$.
Thus, in this case, we obtain
\begin{equation}
p\otimes \fbox{1}\otimes\fbox{111}\otimes\fbox{1111}
\otimes\fbox{1}^{\,\otimes 5}
\end{equation}
as an image of the KKR bijection.

Next, we apply Case 2 procedure to the same modified rigged configuration.
First, we remove the $\mu^{(0)}=(1^{13})$ part from the
quantum space $\mu^{(0)}_+$.
Then we obtain $\fbox{1}^{\,\otimes 13}$ as a part of the image,
and the remaining rigged configuration is
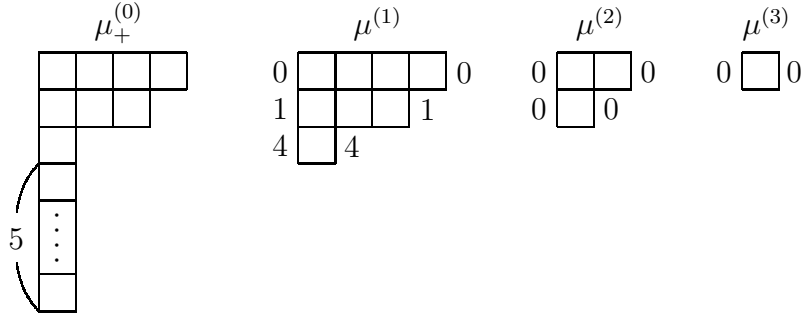
\begin{center}
\unitlength 14pt
\begin{picture}(22,8.5)
\multiput(1,0)(1,0){2}{\line(0,1){4}}
\multiput(1,0)(0,1){2}{\line(1,0){1}}
\multiput(1,3)(0,1){2}{\line(1,0){1}}
\multiput(1.37,1.18)(0,0.4){4}{$\cdot$}
\put(0.17,1.7){5}
\qbezier(1,0)(0.5,0.5)(0.4,1.3)
\qbezier(1,4)(0.5,3.5)(0.4,2.7)
\multiput(1,4)(1,0){2}{\line(0,1){3}}
\put(1,5){\line(1,0){3}}
\multiput(1,6)(0,1){2}{\line(1,0){4}}
\multiput(3,5)(1,0){2}{\line(0,1){2}}
\put(5,6){\line(0,1){1}}
\put(2.5,7.5){$\mu^{(0)}_+$}
\put(8,4){\line(1,0){1}}
\multiput(8,4)(1,0){2}{\line(0,1){3}}
\put(8,5){\line(1,0){3}}
\multiput(8,6)(0,1){2}{\line(1,0){4}}
\multiput(10,5)(1,0){2}{\line(0,1){2}}
\put(12,6){\line(0,1){1}}
\put(7.3,4.2){4}
\put(7.3,5.2){1}
\put(7.3,6.2){0}
\put(9.25,4.17){4}
\put(11.29,5.17){1}
\put(12.28,6.17){0}
\put(9.5,7.5){$\mu^{(1)}$}
\put(15,5){\line(1,0){1}}
\multiput(15,5)(1,0){2}{\line(0,1){2}}
\multiput(15,6)(0,1){2}{\line(1,0){2}}
\put(17,6){\line(0,1){1}}
\put(14.3,5.2){0}
\put(14.3,6.2){0}
\put(16.25,5.2){0}
\put(17.25,6.2){0}
\put(15.5,7.5){$\mu^{(2)}$}
\multiput(20,6)(1,0){2}{\line(0,1){1}}
\multiput(20,6)(0,1){2}{\line(1,0){1}}
\put(19.3,6.2){0}
\put(21.2,6.2){0}
\put(20,7.5){$\mu^{(3)}$}
\end{picture}
\end{center}
At this time, we recognize an implication of the calculation
in Eq.(\ref{555}).
{}From the above diagram we see that the quantum space now is
$\mu^{(0)}_+=\mu^{(1)}\cup (1^{d_1})$,
and all three rows of $\mu^{(1)}$ become simultaneously singular.
This is implied in the following term in Eq.(\ref{555}):
\begin{equation}
d_1=\max\{5,\, 5,\, 5\}.
\end{equation}
(Inside the max term, all factors are 5, and this implies that
all three rows in $\mu^{(1)}$ would simultaneously become singular
when quantum space becomes $\mu^{(0)}_+=\mu^{(1)}\cup (1^{5})$.)

We further proceed along the Case 2 procedure.
As we have said above, we have three possibilities to 
remove a row of $\mu^{(1)}\subset\mu^{(0)}_+$.
Let us remove the row $\{1\}$ from $\mu^{(1)}\subset\mu^{(0)}_+$.
Then we have \fbox{4} as a part of the image, and the remaining
rigged configuration is
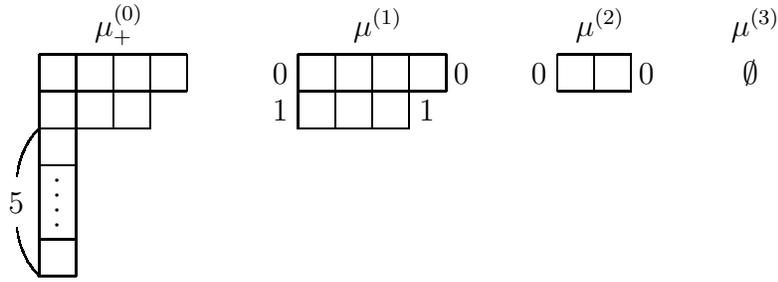
\begin{center}
\unitlength 14pt
\begin{picture}(22,7)
\multiput(1,0)(1,0){2}{\line(0,1){4}}
\multiput(1,0)(0,1){2}{\line(1,0){1}}
\multiput(1,3)(0,1){2}{\line(1,0){1}}
\multiput(1.37,1.18)(0,0.4){4}{$\cdot$}
\put(0.17,1.7){5}
\qbezier(1,0)(0.5,0.5)(0.4,1.3)
\qbezier(1,4)(0.5,3.5)(0.4,2.7)
\put(1,4){\line(1,0){3}}
\multiput(1,4)(1,0){4}{\line(0,1){2}}
\multiput(1,5)(0,1){2}{\line(1,0){4}}
\put(5,5){\line(0,1){1}}
\put(2.5,6.5){$\mu^{(0)}_+$}
\put(8,4){\line(1,0){3}}
\multiput(8,4)(1,0){4}{\line(0,1){2}}
\multiput(8,5)(0,1){2}{\line(1,0){4}}
\put(12,5){\line(0,1){1}}
\put(7.3,4.2){1}
\put(7.3,5.2){0}
\put(11.27,4.2){1}
\put(12.2,5.2){0}
\put(9.5,6.5){$\mu^{(1)}$}
\multiput(15,5)(1,0){3}{\line(0,1){1}}
\multiput(15,5)(0,1){2}{\line(1,0){2}}
\put(14.3,5.2){0}
\put(17.2,5.2){0}
\put(15.5,6.5){$\mu^{(2)}$}
\put(20,5.2){$\emptyset$}
\put(19.7,6.5){$\mu^{(3)}$}
\end{picture}
\end{center}
Again, we encounter two possibilities to remove a row from
$\mu^{(1)}\subset\mu^{(0)}_+$
with now
\begin{equation}
\mu^{(0)}_+=\mu^{(1)}\cup (1^{d_2})=\mu^{(1)}\cup (1^{5}).
\end{equation}
We can infer this by applying Proposition \ref{Q-Q}:
\begin{eqnarray}
d_2&=&\max\{Q^{(1)}_3-Q^{(2)}_3+r^{(1)}_3,\,
Q^{(1)}_4-Q^{(2)}_4+r^{(1)}_4\}
\nonumber\\
&=&\max\{6-2+1,\, 7-2+0\}\nonumber\\
&=&\max\{5,\, 5\}=5.
\end{eqnarray}

Let us remove the row $\{4\}$ from $\mu^{(1)}\subset\mu^{(0)}_+$
and remove the box $\{1\}$ from $(1^{d_2})\subset\mu^{(0)}_+$.
Then we have $\fbox{1}\otimes\fbox{2233}$ as a part of the KKR image,
and the remaining rigged configuration is
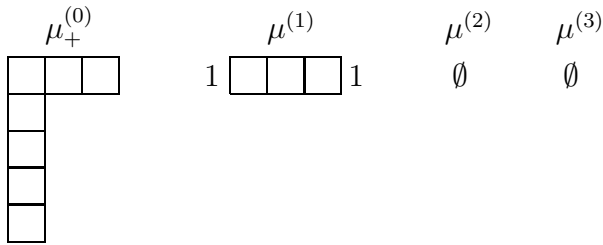
\begin{center}
\unitlength 14pt
\begin{picture}(22,6.5)
\multiput(1,0)(1,0){2}{\line(0,1){5}}
\multiput(1,0)(0,1){4}{\line(1,0){1}}
\multiput(1,4)(0,1){2}{\line(1,0){3}}
\multiput(3,4)(1,0){2}{\line(0,1){1}}
\put(2,5.5){$\mu^{(0)}_+$}
\multiput(7,4)(1,0){4}{\line(0,1){1}}
\multiput(7,4)(0,1){2}{\line(1,0){3}}
\put(6.3,4.2){1}
\put(10.2,4.2){1}
\put(8,5.5){$\mu^{(1)}$}
\put(13,4.2){$\emptyset$}
\put(12.8,5.5){$\mu^{(2)}$}
\put(16,4.2){$\emptyset$}
\put(15.8,5.5){$\mu^{(3)}$}
\end{picture}
\end{center}
At this time, the quantum space $\mu^{(0)}_+$ becomes
$\mu^{(1)}\cup (1^{d_3})$, where we can determine
$d_3$ by Proposition \ref{Q-Q} as follows:
\begin{eqnarray}
d_3&=&\max\{ Q^{(1)}_3-Q^{(2)}_3+r^{(1)}_3\}
\nonumber\\
&=&\max\{ 3-0+1\}=4.
\end{eqnarray}
As a final step of the KKR procedure, 
we remove $\{ 3\}$ from $\mu^{(1)}\subset\mu^{(0)}_+$,
and obtain $\fbox{1}^{\,\otimes 4}\otimes\fbox{222}$
as the rest part of the KKR image.

Both Case 1 and Case 2 procedures above differ from each
other only in the order of removal in the quantum space
$\mu^{(0)}_+$; thus we can apply Kirillov--Schilling--Shimozono's
theorem (Theorem \ref{KSS}) to get the following isomorphism:
\begin{equation}\label{ex:isom}
p\otimes \fbox{1}\otimes\fbox{111}\otimes\fbox{1111}
\otimes\fbox{1}^{\,\otimes 5}
\,\simeq\,
\fbox{1}^{\,\otimes 4}\otimes
\fbox{222}\otimes\fbox{1}\otimes
\fbox{2233}\otimes\fbox{4}\otimes
\left(\fbox{1}^{\,\otimes 13}\right).
\end{equation}
In order to calculate the above isomorphism directly,
we use the diagrammatic method as in (\ref{ex:clclzdon}).
That is, we comapare the right-hand side of Eq.(\ref{ex:isom})
and Eq.(\ref{ex:clclzdon_eq}).
Then we can write down a similar vertex diagram as in
(\ref{ex:clclzdon})
and obtain the left-hand side of (\ref{ex:isom}) as an output.

Below we give a list of all the scattering data obtained
from the rigged configuration here:
\begin{eqnarray*}
s_1&=&\fbox{222}_{\, 4}\otimes\fbox{2233}_{\, 5}\otimes\fbox{4}_{\, 5}\\
s_2&=&\fbox{222}_{\, 4}\otimes\fbox{3}_{\, 5}\otimes\fbox{2234}_{\, 5}\\
s_3&=&\fbox{2222}_{\, 4}\otimes\fbox{233}_{\, 5}\otimes\fbox{4}_{\, 5}\\
s_4&=&\fbox{2222}_{\, 4}\otimes\fbox{3}_{\, 5}\otimes\fbox{234}_{\, 5}
\end{eqnarray*}
$s_1$ is what we have considered in detail.
\rule{5pt}{10pt}
\end{example}
\begin{remark}
By putting letters \fbox{1} on both ends of the above path $p$, we identify
this path as a state of the box-ball systems.
For the sake of simplicity, we tentatively omit frames of tableaux
\fbox{$\ast$} and tensor products $\otimes$, i.e.,
we write the path of the above example as
\begin{equation}
p=1111223214322.
\nonumber
\end{equation}
Then its time evolution is given by
\begin{center}
$t=1$: 1111111122221111332111141111111111111111111111111111111\\*
$t=2$: 1111111111112222111332114111111111111111111111111111111\\*
$t=3$: 1111111111111111222211332411111111111111111111111111111\\*
$t=4$: 1111111111111111111122221343211111111111111111111111111\\*
$t=5$: 1111111111111111111111112232143221111111111111111111111\\*
$t=6$: 1111111111111111111111111121322114322111111111111111111\\*
$t=7$: 1111111111111111111111111112111322111432211111111111111\\*
$t=8$: 1111111111111111111111111111211111322111143221111111111
\end{center}
$t=5$ corresponds to the original path $p$.
We see that there are three solitons of length 4, 3, and 1.
Compare this with the lengths of rows of $\mu^{(1)}$ above.
Furthermore, compare the scattering data $s_3$ at the end of the above
example and $t=1$ path.
Then we see that these three solitons coincide with
three tableaux of $s_3$.
This is the origin of the term ``scattering data.''
In this setting, normal ordering is the way to obtain physically
correct scattering data.
\rule{5pt}{10pt}
\end{remark}

\section{Mode formula and collision states}\label{sec:mode}
In the previous section, we introduce the KKR normal ordered product
\begin{equation}
b_N[d_N]\otimes\cdots\otimes b_1[d_1].
\end{equation}
In order to determine the ordering of sequence $b_1,\cdots ,b_N$ and
associated integers $d_1,\cdots d_N$, we have to refer the
actual KKR procedure.
In this and subsequent sections, we determine these remaining point
by purely crystal bases theoretic scheme.

In the KKR normal ordered product $b_N[d_N]\otimes\cdots\otimes b_1[d_1]$,
the mode $d_1$ is determined by the following simple formula.
This formula also determines the corresponding row $\mu^{(a)}_i$
from which tableau $b_1$ arises.
\begin{proposition}\label{Q-Q}
For the KKR normal ordered product
$b_N[d_N]\otimes\cdots\otimes b_1[d_1]$
we obtain from the
rigged configuration RC${}^{(a-1)}_+$ (Eq.(\ref{levela+})), the mode $d_1$
has the following expression:
\begin{equation}
d_1=\max_{1\leq i\leq l(\mu^{(a)})}\left\{
Q^{(a)}_{\mu^{(a)}_i}-
Q^{(a+1)}_{\mu^{(a)}_i}+
r_i^{(a)}
\right\}\, .
\end{equation}
\end{proposition}
{\bf Proof.}
Consider the KKR bijection on rigged configuration RC${}^{(a-1)}_+$.
We have taken
\begin{equation}
\mu^{(a-1)}_+=\mu^{(a-1)}\cup\mu^{(a)}\cup (1^{d_1}),
\end{equation}
assuming that, while removing $\mu^{(a-1)}$,
the configuration $\mu^{(a)}$ never becomes singular.
Remove $\mu^{(a-1)}$ from the quantum space $\mu^{(a-1)}_+$.
Then we choose $d_1$ so that, just after removing $\mu^{(a-1)}$,
some singular rows appear in $\mu^{(a)}$ for the first time.

We now determine this $d_1$.
To do this, we take arbitrary row $\mu^{(a)}_i$ in
the configuration $\mu^{(a)}$.
Then the condition that this
row becomes singular when we have just removed $\mu^{(a-1)}$ from
$\mu^{(a-1)}_+$ is
\begin{equation}
\left( d_1+Q^{(a)}_{\mu^{(a)}_i}\right)
-2Q^{(a)}_{\mu^{(a)}_i}
+Q^{(a+1)}_{\mu^{(a)}_i}
=r_i^{(a)},
\end{equation}
i.e., the vacancy number of this row
is equal to the corresponding rigging at that time.
Thus, we have
\begin{equation}
d_1=Q^{(a)}_{\mu^{(a)}_i}-
Q^{(a+1)}_{\mu^{(a)}_i}+
r_i^{(a)}.
\end{equation}
These $d_1$'s have different values for different rows $\mu^{(a)}_i$.
Since we define $d_1$ so that the corresponding row is the
first to become singular, we have to take the
maximum of these $d_1$'s.
This completes the proof of the proposition.
\rule{5pt}{10pt}\vspace{4mm}

As a consequence of this formula,
we can derive the following linear dependence of modes $d_i$
on corresponding rigging $r_i^{(a)}$.
\begin{lemma}\label{linear}
Suppose that we have the following KKR normal ordered product
from rigged configuration RC${}^{(a-1)}_+$:
\begin{equation}
b_N[d_N]\otimes\cdots\otimes
b_k[d_k]\otimes\cdots\otimes
b_1[d_1],
\end{equation}
where the tableau $b_k$ originates from the row $\mu^{(a)}_k$,
and the corresponding rigging is $r_k^{(a)}$.
Now we change the rigging $r_k^{(a)}$ to $r_k^{(a)}+1$ and construct
a KKR normal ordered product.
If we can take the ordering of this product to be
$b_N,\cdots ,b_k,\cdots ,b_1$ again,
then the KKR normal ordering is
\begin{equation}\label{eq:lemma_linear}
b_N[d_N]\otimes\cdots\otimes
b_k[d_k+1]\otimes\cdots\otimes
b_1[d_1],
\end{equation}
i.e., $d_j$ $(j\neq k)$ remain the same, and $d_k$ becomes $d_k+1$.
\end{lemma}
{\bf Proof.}
Without loss of generality, we can take $k=2$.
{}From the assumption that we have $b_1[d_1]$ at the right end
of Eq.(\ref{eq:lemma_linear}),
we see that the mode $d_1$ does not change
after we change the rigging $r_2^{(a)}$
(see Proposition \ref{Q-Q}).
We do a KKR procedure in the way described in Case 2 of the previous section.
We first remove the row $\mu^{(a)}_1$, and
next remove the $(1^{d_1})$ part of the quantum space $\mu^{(a-1)}_+$
until some singular rows appear in the configuration $\mu^{(a)}$.
At this time, we can apply Proposition \ref{Q-Q} again to
this removed rigged configuration and obtain the next mode.
Since in formula of Proposition \ref{Q-Q} we take the maximum
of terms, the term corresponding to the row $\mu^{(a)}_2$ is the
maximum one before we change the rigging $r_2^{(a)}$.
Hence it still contributes as the maximal element
even if the rigging becomes $r_2^{(a)}+1$.
{}From this we deduce that the next mode is $d_2+1$.
After removing the row $\mu^{(a)}_2$ and one more box from the quantum space
by the KKR procedure, then the rest of the rigged configuration
goes back to the original situation so that other terms
in the KKR normal ordered product is not different from the original one.
\rule{5pt}{10pt}\vspace{4mm}

To determine modes $d_i$'s,
it is convenient to consider the following state.
\begin{definition}
Consider the KKR bijection on RC${}^{(a-1)}_+$. Remove rows of
the quantum space $\mu^{(a-1)}_+$ by Case 2 procedure
in the previous section, i.e.,
we first remove $\mu^{(a-1)}$ from $\mu^{(a-1)}_+$.
While removing the $(1^d)$ part of the quantum space,
if more than one row of the configuration $\mu^{(a)}$
become simultaneously singular, then we define that
these rows are in collision state.
\rule{5pt}{10pt}
\end{definition}
We choose one of the KKR normal ordered products and fix it.
Suppose that the rightmost elements of it
is $\cdots\otimes B\otimes A$. Then we have the following:
\begin{lemma}
We can always make $B$ and $A$ in collision state by changing
a rigging $r_B$ attached to row $B$.
\end{lemma}
{\bf Proof.}
Let $|A|$ be the width of a tableau $A$.
We can apply the above Lemma \ref{linear}
to make, without changing the other part of the KKR normal ordered product,
\begin{equation}
d_1=
Q^{(a)}_{|A|}-
Q^{(a+1)}_{|A|}+
r_A^{(a)}=
Q^{(a)}_{|B|}-
Q^{(a+1)}_{|B|}+
r_B^{(a)}\, ,
\end{equation}
so that $A$ and $B$ are in collision state.
\rule{5pt}{10pt}
\begin{example}
Consider $s_1$ and $s_3$ in Example \ref{ex1}.
In $s_1$, \fbox{2233} and \fbox{4} are in the collision state,
and, in $s_3$, \fbox{233} and \fbox{4} are in the collision state.
\rule{5pt}{10pt}
\end{example}

\section{Energy functions and the KKR bijection}\label{sec:energy&KKR}
In the previous sections, we give crystal interpretation for several
properties of the KKR bijection, especially with respect to combinatorial
$R$ matrices.
Now it is a point to determine all modes $d_i$
in scattering data by use of the $H$ function or the energy function 
of a product $B\otimes A$.
We consider the rigged configuration RC${}^{(a-1)}_+$
(Eq.(\ref{levela+}));
that is, its quantum space is
\begin{equation}
\mu^{(a-1)}_+=\mu^{(a-1)}\cup\mu^{(a)}\cup (1^{d_1}).
\end{equation}
In the following discussion, we take $a=1$ without loss of generality
and remove $\mu^{(0)}$ as a first step.

To describe the main result, we prepare some conventions and notation.
For the KKR normal ordered product, we denote the rightmost part
as $\cdots B[d_2]\otimes A[d_1]$,
where the lengths of tableaux are $|A|=L$ and $|B|=M$.
Tableaux $A$ and $B$ originate from rows of $\mu^{(1)}$,
which we also denote as row $A$ and row $B$ for the sake of simplicity.
The difference of $Q^{(i)}_j$'s before and after the removal of row $A$
is $\Delta Q^{(i)}_j$, i.e.,
\begin{equation}
\Delta Q^{(i)}_j:=(Q^{(i)}_j\mbox{ just before removal of $A$})
-(Q^{(i)}_j\mbox{ just after removal of $A$}).
\end{equation}
Then we have the following theorem.
\begin{theorem}\label{unwinding}
If $A$ and $B$ in the KKR normal ordered product are successive
(i.e. $\cdots B[d_2]\otimes A[d_1]$), then we have
\begin{equation}
\Delta Q^{(2)}_M=\mbox{the unwinding number of }B\otimes A.
\end{equation}
\end{theorem}
{\bf Proof} will be given in the next section.
\rule{5pt}{10pt}\vspace{4mm}\\
We give two examples of this theorem.
\begin{example}\label{ex:unwinding1}
Consider the following $\mathfrak{sl}_5$ rigged configuration:
\begin{equation}
\mu^{(0)}_+=\{1,1,1^3,3,3,1,3,1,4\}
\nonumber
\end{equation}
\begin{center}
\unitlength 14pt
\begin{picture}(22,6)
\put(0,0){\line(1,0){1}}
\multiput(0,0)(1,0){2}{\line(0,1){5}}
\multiput(0,1)(0,1){3}{\line(1,0){3}}
\multiput(2,1)(1,0){2}{\line(0,1){4}}
\multiput(0,4)(0,1){2}{\line(1,0){4}}
\put(4,4){\line(0,1){1}}
\put(-0.7,0.2){4}
\multiput(-0.7,1.2)(0,1){4}{0}
\put(1.2,0.17){0}
\multiput(3.2,1.17)(0,1){3}{0}
\put(4.2,4.17){0}
\put(1.5,5.5){$\mu^{(1)}$}
\put(7,2){\line(1,0){1}}
\multiput(7,2)(1,0){2}{\line(0,1){3}}
\put(7,3){\line(1,0){3}}
\multiput(7,4)(0,1){2}{\line(1,0){4}}
\multiput(9,3)(1,0){2}{\line(0,1){2}}
\put(11,4){\line(0,1){1}}
\put(6.3,2.2){1}
\put(6.3,3.2){2}
\put(6.3,4.2){1}
\put(8.23,2.17){1}
\put(10.23,3.17){1}
\put(11.23,4.17){0}
\multiput(9.18,4.25)(1,0){2}{$\times$}
\put(7.18,2.25){$\times$}
\put(8.5,5.5){$\mu^{(2)}$}
\multiput(7,0)(0,0.4){15}{\line(0,1){0.23}}
\multiput(10,0)(0,0.4){15}{\line(0,1){0.23}}
\put(7,0.9){\line(1,0){1}}
\put(9,0.9){\line(1,0){1}}
\put(6.89,0.665){$<$}
\put(9.44,0.665){$>$}
\put(8.26,0.6){3}
\put(14,3){\line(1,0){1}}
\multiput(14,3)(1,0){2}{\line(0,1){2}}
\multiput(14,4)(0,1){2}{\line(1,0){2}}
\put(16,4){\line(0,1){1}}
\put(13.3,3.2){0}
\put(13.3,4.2){0}
\put(15.2,3.17){0}
\put(16.2,4.17){0}
\put(14.5,5.5){$\mu^{(3)}$}
\multiput(19,4)(1,0){2}{\line(0,1){1}}
\multiput(19,4)(0,1){2}{\line(1,0){1}}
\put(18.3,4.2){0}
\put(20.25,4.2){0}
\put(19,5.5){$\mu^{(4)}$}
\end{picture}
\end{center}
We have assumed that we had already removed the $\mu^{(0)}$
part from $\mu^{(0)}_+$ (the expression of $\mu^{(0)}_+$
is reordered form, see Step 1 of Definition \ref{defKKR}).
Then, by the KKR procedure, we remove rows of $\mu^{(0)}_+$
from right to left in the above ordering
and obtain the following KKR normal ordered product:
\begin{equation}
\fbox{1}\otimes\fbox{2}\otimes\fbox{1}^{\,\otimes 3}\otimes
\fbox{222}\otimes\fbox{333}\otimes\fbox{1}\otimes\fbox{244}
\otimes\fbox{1}\otimes\fbox{2335}\, .
\end{equation}
The rightmost part of the product satisfies
\begin{equation}
\mbox{the unwinding number of }\,\fbox{244}\otimes\fbox{2335}=2.
\end{equation}

In the above diagram, boxes with cross
\unitlength 13pt
\begin{picture}(1,1)(0,0.2)
\multiput(0,0)(1,0){2}{\line(0,1){1}}
\multiput(0,0)(0,1){2}{\line(1,0){1}}
\put(0.15,0.25){$\times$}
\end{picture} in $\mu^{(2)}$
mean that when we obtain \fbox{2335} these boxes are removed
by the KKR procedure.
Since the width of \fbox{244} is 3, we have $\Delta Q^{(2)}_3=2$,
which agrees with Theorem \ref{unwinding}.
\rule{5pt}{10pt}
\end{example}
\begin{example}\label{ex:unwinding2}
Let us consider one more elaborated example in $\mathfrak{sl}_6$.
\begin{equation}
\mu^{(0)}_+=\{1^4,4,1^3,4,1^3,6,1,6,1^4,8,1^2,8\}
\nonumber
\end{equation}
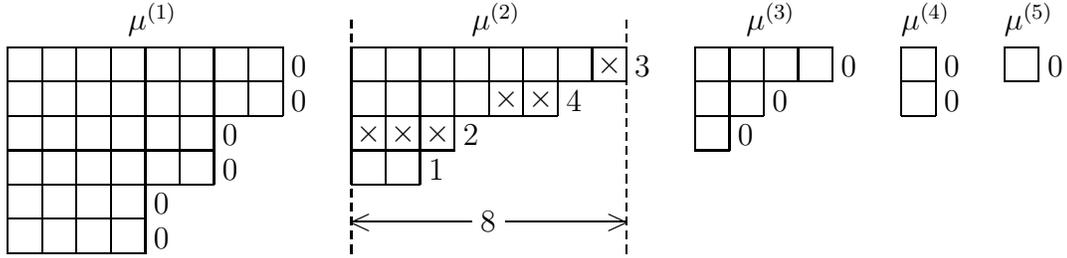
\begin{center}
\unitlength 13pt
\begin{picture}(32,7)
\multiput(1,0)(1,0){5}{\line(0,1){6}}
\multiput(1,0)(0,1){2}{\line(1,0){4}}
\multiput(1,2)(0,1){2}{\line(1,0){6}}
\multiput(1,4)(0,1){3}{\line(1,0){8}}
\multiput(6,2)(1,0){2}{\line(0,1){4}}
\multiput(8,4)(1,0){2}{\line(0,1){2}}
\multiput(5.25,0.15)(0,1){2}{0}
\multiput(7.25,2.15)(0,1){2}{0}
\multiput(9.25,4.15)(0,1){2}{0}
\put(4.5,6.5){$\mu^{(1)}$}
\put(11,2){\line(1,0){2}}
\multiput(11,2)(1,0){3}{\line(0,1){4}}
\put(11,3){\line(1,0){3}}
\put(11,4){\line(1,0){6}}
\multiput(11,5)(0,1){2}{\line(1,0){8}}
\put(14,3){\line(0,1){3}}
\multiput(15,4)(1,0){3}{\line(0,1){2}}
\multiput(18,5)(1,0){2}{\line(0,1){1}}
\put(13.25,2.15){1}
\put(14.25,3.15){2}
\put(17.25,4.15){4}
\put(19.25,5.15){3}
\put(14.5,6.5){$\mu^{(2)}$}
\multiput(11.15,3.25)(1,0){3}{$\times$}
\multiput(15.15,4.25)(1,0){2}{$\times$}
\put(18.15,5.25){$\times$}
\multiput(11,0)(0,0.41){17}{\line(0,1){0.24}}
\multiput(19,0)(0,0.41){17}{\line(0,1){0.24}}
\put(11,0.93){\line(1,0){3.5}}
\put(15.5,0.93){\line(1,0){3.5}}
\put(14.75,0.63){8}
\put(10.9,0.675){$<$}
\put(18.4,0.675){$>$}
\put(21,3){\line(1,0){1}}
\multiput(21,3)(1,0){2}{\line(0,1){3}}
\put(21,4){\line(1,0){2}}
\put(23,4){\line(0,1){2}}
\multiput(21,5)(0,1){2}{\line(1,0){4}}
\multiput(24,5)(1,0){2}{\line(0,1){1}}
\put(22.25,3.15){0}
\put(23.25,4.15){0}
\put(25.25,5.15){0}
\put(22.5,6.5){$\mu^{(3)}$}
\multiput(27,4)(1,0){2}{\line(0,1){2}}
\multiput(27,4)(0,1){3}{\line(1,0){1}}
\multiput(28.25,4.15)(0,1){2}{0}
\put(27,6.5){$\mu^{(4)}$}
\multiput(30,5)(1,0){2}{\line(0,1){1}}
\multiput(30,5)(0,1){2}{\line(1,0){1}}
\put(31.25,5.15){0}
\put(30,6.5){$\mu^{(5)}$}
\end{picture}
\end{center}
We have suppressed the vacancy numbers for the sake
of simplicity.
By the KKR procedure, we have the following
KKR normal ordered product:
\begin{eqnarray}
&&
\fbox{1}^{\,\otimes 4}\otimes\fbox{2222}\otimes
\fbox{1}^{\,\otimes 3}\otimes\fbox{2223}\otimes
\fbox{1}^{\,\otimes 3}\otimes\fbox{222334}\otimes
\fbox{1}\otimes\fbox{233344}\nonumber\\
&&\hspace{3mm}\otimes\,
\fbox{1}^{\,\otimes 4}\otimes\fbox{22223345}\otimes
\fbox{1}^{\,\otimes 2}\otimes\fbox{22333346}\, .
\end{eqnarray}
The rightmost part of this product satisfies
\begin{equation}
\mbox{the unwinding number of }\,
\fbox{22223345}\otimes\fbox{22333346}=6.
\end{equation}
Since the width of \fbox{22223345} is 8,
this means that $\Delta Q^{(2)}_8=6$,
and this agrees with the number of
\unitlength 13pt
\begin{picture}(1,1)(0,0.2)
\multiput(0,0)(1,0){2}{\line(0,1){1}}
\multiput(0,0)(0,1){2}{\line(1,0){1}}
\put(0.15,0.25){$\times$}
\end{picture}
in $\mu^{(2)}|_{\leq 8}$ of the above diagram.
\rule{5pt}{10pt}
\end{example}

Implication of this theorem is as follows.
Without loss of generality,
we take $A$ and $B$ as a collision state.
We are choosing a normalization for the $H$ function as
\begin{equation}
H:=\mbox{the winding number of }B\otimes A.
\end{equation}
By the above definition, $\Delta Q^{(i)}_M$ is equal to the number
of boxes which are removed from $\mu^{(i)}|_{\leq M}$ when we
remove the row $A$ by the KKR procedure.
Thus, if we remove row $A\subset \mu^{(1)}$, then $\Delta Q^{(1)}_M$ is
(recall that $M:=|B|$ and $L:=|A|$)
\begin{equation}
\Delta Q^{(1)}_M =\min \{ M,L\} .
\end{equation}
{}From the above theorem we have
\begin{eqnarray}
\Delta Q^{(1)}_M-\Delta Q^{(2)}_M&=&
  \min \{ M,L\}-\mbox{the unwinding number of }B\otimes A\nonumber\\
&=&\mbox{the winding number of }B\otimes A\nonumber\\
&=&H(B\otimes A).
\end{eqnarray}

On the other hand, since $A$ and $B$ are in the collision state, from
Proposition \ref{Q-Q} we have
\begin{equation}
d_1=
Q^{(1)}_{|A|}-Q^{(2)}_{|A|}+r_A=
Q^{(1)}_{|B|}-Q^{(2)}_{|B|}+r_B\,
\end{equation}
just before we remove the row $A$.
After removing the row $A$, we again apply Proposition \ref{Q-Q}
to the rigged configuration which has been modified by removal of the row $A$
under the KKR procedure.
We then have
\begin{equation}
d_2=Q^{(1)}_{|B|}-Q^{(2)}_{|B|}+r_B
\end{equation}
after removing $A$.
By the definition of $\Delta Q^{(i)}_j$, we have
\begin{equation}
d_1-d_2=\Delta Q^{(1)}_M-\Delta Q^{(2)}_M\, ,
\end{equation}
so that, combining the above arguments, we obtain
\begin{equation}
d_1-d_2=H\, ,
\end{equation}
as a consequence of the above theorem.
\begin{remark}
{}From the point of view of the box-ball systems,
Theorem \ref{unwinding} thus asserts that
minimal separation between two successive solitons
is equal to the energy function of the both solitons.
\rule{5pt}{10pt}
\end{remark}

\noindent
{\bf Proof of Theorem \ref{clcl-zdon}.}
What we have to do is to identify the KKR normal ordered product
(Definition \ref{KKRnorder})
with normal ordering (Definition \ref{norder})
which is defined in terms of the crystal base theory.
Again, we consider the rigged configuration RC${}^{(a-1)}_+$
(Eq.(\ref{levela+}))
with $a=1$.

First, we give an interpretation of the above arguments
about the collision states in terms of
the normal ordering of Definition \ref{norder}.
Consider the following normal ordered product in the sense of
Definition \ref{norder}:
\begin{equation}
b_N[d_N]\otimes\cdots\otimes
b_2[d_2]\otimes b_1[d_1].
\end{equation}
Concentrate on a particular successive pair
$b_{i+1}[d_{i+1}]\otimes b_i[d_i]$ within this scattering data.
The isomorphism of the affine combinatorial $R$ then gives
\begin{equation}
b_{i+1}[d_{i+1}]\otimes b_i[d_i]
\,\simeq\,
b_i'[d_i-H]\otimes b_{i+1}'[d_{i+1}+H],
\end{equation}
where $H$ is a value of $H$ function on this product,
and $b_{i+1}\otimes b_i\simeq b_i'\otimes b_{i+1}'$
is the corresponding isomorphism under the
classical combinatorial $R$ matrix.
Since the modes $d_i$ depend linearly on the corresponding
rigging $r_i$ (see Eq.(\ref{defmode})),
we can adjust $r_i$ to make that both
\begin{equation}
\cdots\otimes
b_{i+1}[d_{i+1}]\otimes
b_i[d_i]\otimes\cdots\quad
{\rm and}\quad
\cdots\otimes
b_i'[d_i-H]\otimes b_{i+1}'[d_{i+1}+H]\otimes\cdots
\end{equation}
are simultaneously normal ordered,
where the abbreviated parts in the above expression are unchanged.
{}From Definition \ref{norder} we see that all
normal ordered products possess the common set of modes $\{ d_i\}$.
Thus, if this adjustment is already taken into account,
then the modes $d_{i+1}$ and $d_i$ satisfy
\begin{equation}
d_i-d_{i+1}=H,
\end{equation}
which is the same relation as what we have seen in the case
of KKR normal ordering.

To summarize, both the KKR normal ordering and
normal ordering share the following common properties:
\begin{enumerate}
\item Each $b_i$ is a tableau which is obtained as
 a KKR image of the rigged configuration RC${}^{(a)}$ (Eq.(\ref{levela}))
 with $a=1$.
 They commute with each other under the isomorphism of
 $\mathfrak{sl}_{n-1}$ combinatorial
 $R$ matrices with letters from 2 to $n$.
\item Consider a normal ordered product.
 If we can change some riggings $r_i$ without changing the order
 of elements in normal ordering, then
 each mode $d_i$ depends linearly on the corresponding
 rigging $r_i$.
\item Concentrate on a particular product
 $b_i\otimes b_j$
 inside a scattering data, then we can adjust corresponding
 rigging $r_i$ to make that both
 \begin{equation}
 \cdots\otimes
 b_i[d_i]\otimes b_j[d_j]
 \otimes\cdots\, 
 \quad\mbox{ and }\quad
 \cdots\otimes
 b_j'[d_j']\otimes b_i'[d_i']
 \otimes\cdots\, ,\nonumber
 \end{equation}
 where
 $b_i\otimes b_j
 \stackrel{R}{\simeq}
 b_{j}'\otimes b_{i}'$,
 are simultaneously normal ordered.
 If we have already adjusted the rigging $r_i$
 in such a way, then the difference of the successive modes $d_i$ and $d_j$
 is equal to
 \begin{equation}
 d_j-d_i=H\, ,
 \end{equation}
 i.e., the value of the $H$ function on this product.
\end{enumerate}

{}From these observations we see that the both modes $d_i$ defined by
Proposition \ref{Q-Q} and Eq.(\ref{defmode})
are essentially identical.
Thus the KKR normal ordered products are normal ordered products
in the sense of Definition \ref{norder}.
On the contrary, we can say that all the normal ordered products
are, in fact, KKR normal ordered.
To see this, take one of the normal ordered products
\begin{equation}
b_N[d_N]\otimes\cdots\otimes
b_2[d_2]\otimes b_1[d_1]\in \mathcal{S}_1,
\end{equation}
where $\mathcal{S}_1$ is defined in Definition \ref{norder}.
{}From this scattering data we construct the element
\begin{equation}
\fbox{1}^{\,\otimes d_N}\otimes b_N\otimes
\fbox{1}^{\,\otimes (d_{N-1}-d_N)}\otimes\cdots\otimes
b_2\otimes
\fbox{1}^{\,\otimes (d_1-d_2)}\otimes b_1.
\end{equation}
Then, in view of the isomorphism of affine combinatorial $R$ matrices,
each power $d_{i-1}-d_i$
is larger than the corresponding $H$ function
(because if it is not the case, then we can permute
$b_i[d_i]\otimes b_{i-1}[d_{i-1}]$ to make the $(i-1)$th mode as
$d_i +H$, i.e., larger than the original $d_{i-1}$
in contradiction to the definition of the normal ordering).
{}From Theorem \ref{unwinding} and the
arguments following it we see that it is a
sufficient condition to be a KKR normal ordered product
(by a suitable choice of riggings $r^{(1)}_i$; since other
information, i.e., RC${}^{(1)}$ can be
determined from $b_N, \cdots, b_1$ in the given scattering data).
Thus we can apply the inverse of the KKR bijection and obtain the
corresponding rigged configuration.

In the earlier arguments, we have interpreted the $\Phi^{(a)}$
operator in terms of the KKR bijection (Proposition \ref{moritsuke1}).
Now we interpret the $C^{(a)}$ operators
or, in other words, the normal ordering in terms of the KKR bijection.

Hence we complete the proof of Theorem \ref{clcl-zdon}.
\rule{5pt}{10pt}

\begin{remark}\label{characterization}
In the above arguments, we see that the normal ordered
scattering data can be identified with the paths
obtained from the rigged configuration RC${}^{(a-1)}_+$.
In particular, if the element
\begin{equation}
s=b_N[d_N]\otimes\cdots\otimes b_2[d_2]\otimes
b_1[d_1]
\end{equation}
satisfies the two conditions
\begin{enumerate}
\item $b_N\otimes\cdots\otimes b_2\otimes b_1$ is a path of
      RC${}^{(a)}$,
\item every difference of modes $d_i$ satisfies
      the condition $d_{i}-d_{i+1}\geq H(b_{i+1}\otimes b_i)$,
\end{enumerate}
then $s$ can be realized as an image of RC${}^{(a-1)}_+$.
Therefore we obtain the following characterization of
normal orderings.

Let $\mathcal{S}_{N+1}$ be the set defined in Definition \ref{norder}.
Consider an element
$s=b_N[d_N]\otimes\cdots\otimes b_2[d_2]\otimes
b_1[d_1]\in\mathcal{S}_{N+1}$.
Then $s\in\mathcal{S}_1$ if and only if
the modes $d_i$ of $s$ satisfy
the condition $d_{i}-d_{i+1}\geq H(b_{i+1}\otimes b_i)$
for all $1\leq i\leq N-1$.
\rule{5pt}{10pt}
\end{remark}

\section{Proof of Theorem \ref{unwinding}}\label{sec:proof}
Proof of the theorem is divided into 6 steps.\vspace{4mm}

\noindent
{\bf Step 1:} Let us introduce some notation used throughout the proof.
Consider the rigged configuration RC${}^{(a-1)}_+$ (Eq.(\ref{levela+}))
with $a=1$.
Let the rightmost elements of a KKR normal ordered product be
$\cdots\otimes B\otimes A$.
When we remove the $k$th box from the right end of the row $A\subset\mu^{(1)}$,
we remove the box $\alpha^{(j)}_k$ from the configuration $\mu^{(j)}$
by the KKR bijection.
That is, when we remove the $k$th box of a row $A$, the boxes
\begin{equation}
\alpha_k^{(2)},\,\alpha_k^{(3)},\,
\cdots ,\,\alpha_k^{(n-1)}
\end{equation}
are also removed.
In some cases, we have
\begin{equation}
\alpha_k^{(j-1)}\neq\emptyset ,\,
\alpha_k^{(j)}=\emptyset ,\,
\alpha_k^{(j+1)}=\emptyset ,\,
\cdots ,
\end{equation}
for some $j\leq n-1$.
The box adjacent to the left of the box $\alpha^{(j)}_k$
is the box $\alpha^{(j)}_k-1$.
We sometimes express a row by its rightmost box.
Then we have the following:
\begin{lemma}\label{tanchousei}
For fixed $j$, $col(\alpha^{(j)}_k)$ monotonously decrease
with respect to $k$, i.e., $col(\alpha^{(j)}_k)>col(\alpha^{(j)}_{k+1})$.
\end{lemma}
{\bf Proof.}
First, we consider $col(\alpha^{(2)}_k)$.
When we remove $\alpha^{(1)}_k$, i.e.,
the $k$th box from right end of a row $A\subset\mu^{(1)}$,
then we remove the box $\alpha^{(2)}_k$ from $\mu^{(2)}$
and continue as far as possible.
In the next step, we remove the box $\alpha^{(1)}_{k+1}$
from the row $A$ which satisfies
\begin{equation}
col(\alpha^{(1)}_{k+1})=
col(\alpha^{(1)}_{k})-1.
\end{equation}
After the box $\alpha^{(1)}_{k+1}$, we remove a box
$\alpha^{(2)}_{k+1}$, which has the following two possibilities:
\begin{enumerate}
\item[(1)] $\alpha^{(2)}_{k+1}$ and $\alpha^{(2)}_{k}$
are on the same row, or
\item[(2)] $\alpha^{(2)}_{k+1}$ and $\alpha^{(2)}_{k}$
are on different rows.
\end{enumerate}
In case (1), we have 
$col(\alpha^{(2)}_{k+1})=
col(\alpha^{(2)}_{k})-1$.
In case (2), we have
$col(\alpha^{(2)}_{k+1})<
col(\alpha^{(2)}_{k})-1$,
since if $col(\alpha^{(2)}_{k+1})=
col(\alpha^{(2)}_{k})-1$,
then we can choose $\alpha^{(2)}_{k+1}$ from the same
row with $\alpha^{(2)}_{k}$.
In both cases, $col(\alpha^{(2)}_k)$ monotonously decreases
with respect to $k$.

In the same way, we assume that until some $j$, $col(\alpha^{(j)}_k)$
monotonously decreases with respect to $k$.
Then, from the relation
\begin{equation}
col(\alpha^{(j)}_{k+1})\leq
col(\alpha^{(j)}_k)-1,
\end{equation}
we can show that $col(\alpha^{(j+1)}_k)$ also monotonously decreases
with respect to $k$.
By induction, this gives a proof of the lemma.
\rule{5pt}{10pt}\vspace{4mm}

\noindent
{\bf Step 2:}
When we remove boxes $\alpha^{(1)}_k,$ $\alpha^{(2)}_k,$
$\alpha^{(3)}_k,\cdots$,
the vacancy numbers of the rigged configuration change
in a specific way.
In this step, we pursue this characteristic pattern
before and after the removal.

\begin{figure}
\begin{center}
\unitlength 10pt
\begin{picture}(30,12)
\multiput(-4,5)(0.5,0){4}{$\cdot$}
\put(0,0){\line(0,1){10}}
\put(0,10){\line(1,0){7}}
\qbezier(0,0)(5,0.5)(7,10)
\put(2,11){$\mu^{(i-1)}$}
\put(2.5,1.5){$\times$}
\put(1,2.6){$\alpha^{(i-1)}_k$}
\multiput(3.3,1.8)(0.41,0){28}{\line(1,0){0.23}}
\thicklines
\qbezier(4.5,1.8)(7.5,1.8)(8.5,5)
\qbezier(8.5,5)(9.5,7.95)(12,7.95)
\qbezier(12,7.95)(15,7.95)(20,7.95)
\thinlines
\put(12,0){\line(0,1){10}}
\put(12,10){\line(1,0){7}}
\qbezier(12,0)(17,0.5)(19,10)
\put(14.5,11){$\mu^{(i)}$}
\put(14.5,1.5){$\times$}
\put(13.4,2.6){$\alpha^{(i)}_k$}
\multiput(12,7.8)(0.407,0){44}{\line(1,0){0.23}}
\multiput(12,1.9)(0,0.16){37}{\line(1,0){0.25}}
\put(10.2,4.3){+1}
\put(24,0){\line(0,1){10}}
\put(24,10){\line(1,0){7}}
\qbezier(24,0)(29,0.5)(31,10)
\put(26,11){$\mu^{(i+1)}$}
\put(29.6,7.5){$\times$}
\put(27.5,8.5){$\alpha^{(i+1)}_k$}
\multiput(32,5)(0.5,0){4}{$\cdot$}
\end{picture}
\end{center}
\caption{Schematic diagram of Eq.(\ref{aida}).
To make the situation transparent, three Young diagrams
are taken to be the same.
The shaded region of $\mu^{(i)}$ is Eq.(\ref{aida})
whose coquantum numbers are increased by 1.
Imagine that we are removing boxes
$\beta^{(1)}_k$, $\beta^{(2)}_k$, $\beta^{(3)}_k$, $\cdots$
from the left configuration to the right
one according to the KKR procedure (see Step 3).
We can think of it as some kind of a ``particle''
traveling from left to right until stopped.
Then, when we remove a row $B$, the curved thick line in the
figure looks like a ``potential wall'' which prevents the particle from
going rightwards.}
\label{fig1}
\end{figure}
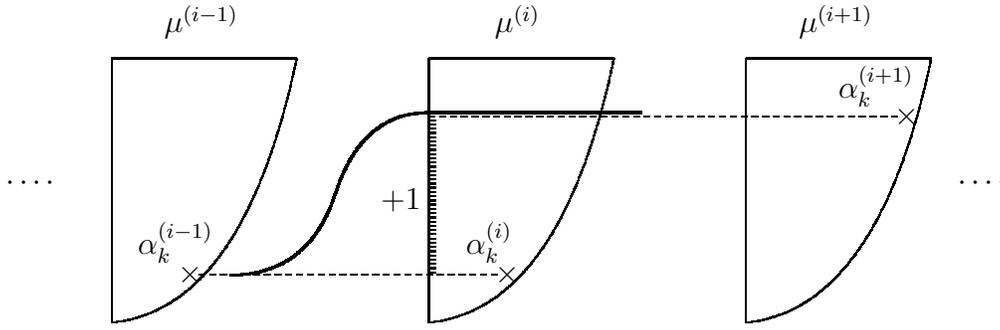

First, consider the case $\alpha^{(i+1)}_k\neq\emptyset$,
i.e., remove the boxes $\alpha^{(1)}_k$, $\alpha^{(2)}_k$,
$\cdots$, $\alpha^{(i+1)}_k$, $\cdots$.
If $col(\alpha^{(i)}_{k})<col(\alpha^{(i+1)}_{k})$,
then the vacancy numbers attached to the rows 
$\alpha^{(i)}(\neq\alpha^{(i)}_{k})$ of
the configuration $\mu^{(i)}$ within the region
\begin{equation}\label{aida}
col(\alpha^{(i)}_{k})\leq col(\alpha^{(i)})<
col(\alpha^{(i+1)}_{k}),
\end{equation}
increase by 1 (see Fig. \ref{fig1}).
To see this, let us tentatively write $col(\alpha^{(i)})=l$.
Recall that the vacancy number $p^{(i)}_l$ for this row is
\begin{equation}
p^{(i)}_l=Q^{(i-1)}_l-2Q^{(i)}_l+Q^{(i+1)}_l.
\end{equation}
After removing boxes $\alpha^{(i-1)}_{k}$, $\alpha^{(i)}_{k}$,
$\alpha^{(i+1)}_{k}$, we see that
$Q^{(i-1)}_l$ and $Q^{(i)}_l$ decrease by 1;
on the other hand, $Q^{(i+1)}_l$ do not change
(because of Eq.(\ref{aida}) combined with
$col(\alpha^{(i-1)}_k)\leq col(\alpha^{(i)}_k)$).
Summing up these contributions, the vacancy numbers
$p^{(i)}_l$ increase by 1.
It also implies that the coquantum numbers (i.e., the vacancy numbers
minus riggings for the corresponding rows) also increase by 1.
Similarly, if we have the condition
$col(\alpha^{(i)}_{k})=col(\alpha^{(i+1)}_{k})$,
then the vacancy numbers for rows $\alpha^{(i)}$ of $\mu^{(i)}$
in the region
\begin{equation}
col(\alpha^{(i)}_{k})\leq col(\alpha^{(i)})
\end{equation}
do not change, since $Q^{(i+1)}_l$ also decrease by 1 in this case.

\begin{figure}
\begin{center}
\unitlength 10pt
\begin{picture}(30,12)
\multiput(-4,5)(0.5,0){4}{$\cdot$}
\put(0,0){\line(0,1){10}}
\put(0,10){\line(1,0){7}}
\qbezier(0,0)(5,0.5)(7,10)
\put(2,11){$\mu^{(i-1)}$}
\put(4,3.5){$\times$}
\multiput(4.75,3.8)(0.4,0){29}{\line(1,0){0.23}}
\put(2.5,4.7){$\alpha^{(i-1)}_k$}
\thicklines
\qbezier(5.3,3.9)(8,3.9)(9,6)
\qbezier(9,6)(9.8,8)(10,10.5)
\thinlines
\put(12,0){\line(0,1){10}}
\put(12,10){\line(1,0){7}}
\qbezier(12,0)(17,0.5)(19,10)
\put(14.5,11){$\mu^{(i)}$}
\put(16,3.5){$\times$}
\put(15,4.7){$\alpha^{(i)}_k$}
\multiput(12,3.9)(0,0.16){38}{\line(1,0){0.25}}
\put(10,6){+1}
\put(24,0){\line(0,1){10}}
\put(24,10){\line(1,0){7}}
\qbezier(24,0)(29,0.5)(31,10)
\put(26,11){$\mu^{(i+1)}$}
\multiput(32,5)(0.5,0){4}{$\cdot$}
\end{picture}
\end{center}
\caption{Schematic diagram of Eq.(\ref{haji}).
Coquantum numbers of the shaded region in $\mu^{(i)}$
are increased by 1 (Eq.(\ref{haji})).
Thick line in the figure shows a ``potential wall''
as in Fig. \ref{fig1}.}
\label{fig2}
\end{figure}
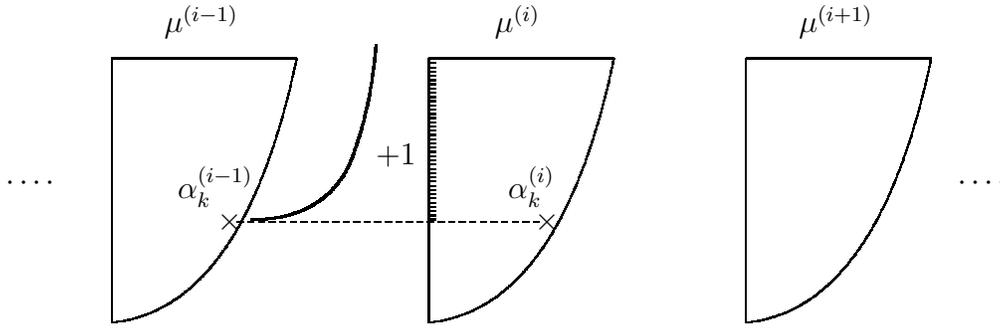

Next, consider the case $\alpha^{(i)}_{k}\neq\emptyset$ and
$\alpha^{(i+1)}_{k}=\emptyset$ with $i\leq n-1$.
Then the vacancy numbers for rows $\alpha^{(i)}(\neq\alpha^{(i)}_{k})$
of $\mu^{(i)}$ in the region
\begin{equation}\label{haji}
col(\alpha^{(i)}_{k})\leq col(\alpha^{(i)})
\end{equation}
increase by 1, since within the vacancy number
$p^{(i)}_l=Q^{(i-1)}_l-2Q^{(i)}_l+Q^{(i+1)}_l$ for a row $\alpha^{(i)}$
(with width $l$),
$Q^{(i-1)}_l$ and $Q^{(i)}_l$ decrease by 1; on the other hand,
$Q^{(i+1)}_l$ do not change, since now $\alpha^{(i+1)}_{k}=\emptyset$
(see Fig.\ref{fig2}).
Therefore its coquantum number also increase by 1.
The above arguments in Step 2 are summarized in (I), (II), and (III) of
Lemma \ref{potential} below.
\vspace{4mm}

\noindent
In the rest of this Step 2,
we show that once regions Eq.(\ref{aida}) or Eq.(\ref{haji}) of $\mu^{(i)}$
become nonsingular in the way described above, then they
never become singular even when we are removing the rest of a row $A$.
To begin with, consider the effect of $\alpha^{(i-1)}_{k+1}$,
$\alpha^{(i)}_{k+1}$, and $\alpha^{(i+1)}_{k+1}$.
In what follows, we first treat $\alpha^{(i+1)}_k\neq\emptyset$,
and then $\alpha^{(i+1)}_k=\emptyset$.

We see that if $\alpha^{(i+1)}_{k}\neq\emptyset$, then
$\alpha^{(i+1)}_{k+1}\neq\emptyset$.
This is because: (1) the row $\alpha^{(i+1)}_{k}-1$ becomes singular,
since we have removed a box $\alpha^{(i+1)}_{k}$, and
(2) the width of the row $\alpha^{(i+1)}_{k}-1$ satisfies
\begin{equation}
col(\alpha^{(i)}_{k+1})\leq col(\alpha^{(i+1)}_{k})-1
=col(\alpha^{(i+1)}_{k}-1)
\end{equation}
because of the relation
\begin{equation}\label{kantan1}
col(\alpha^{(i)}_{k+1})<
col(\alpha^{(i)}_{k})\leq col(\alpha^{(i+1)}_{k})
\end{equation}
(the first $<$ is by Lemma \ref{tanchousei}, and the next $\leq$ is
by the definition of the KKR bijection).
Thus we can remove the row $\alpha^{(i+1)}_{k}-1$ of $\mu^{(i+1)}$
(or, if exists, s shorter singular row) as the next box
$\alpha^{(i+1)}_{k+1}(\neq\emptyset $, as requested).

Now we are assuming that $\alpha^{(i+1)}_{k}\neq\emptyset$ and
thus $\alpha^{(i+1)}_{k+1}\neq\emptyset$.
We have two inequalities
\begin{enumerate}
\item[(i)] $col(\alpha^{(i-1)}_{k+1})<col(\alpha^{(i-1)}_{k})
\leq col(\alpha^{(i)}_{k})$ (by the same reason as Eq.(\ref{kantan1}))
and 
\item[(ii)] $col(\alpha^{(i)}_{k+1})<col(\alpha^{(i)}_{k})$
(by Lemma \ref{tanchousei}).
\end{enumerate}
Using these two relations,
consider the change of the vacancy number corresponding to
the rows $\alpha^{(i)}\in\mu^{(i)}$ within the region
\begin{equation}\label{kantan2}
col(\alpha^{(i)}_{k})\leq col(\alpha^{(i)})<
col(\alpha^{(i+1)}_{k}),
\end{equation}
when we remove boxes $\alpha^{(i-1)}_{k+1}$,
$\alpha^{(i)}_{k+1}$, and $\alpha^{(i+1)}_{k+1}$
(see Fig. \ref{fig1}).
Let us write $col(\alpha^{(i)})=l$, then the vacancy number $p^{(i)}_l$
for a row $\alpha^{(i)}$ is
$p^{(i)}_l=Q^{(i-1)}_l-2Q^{(i)}_l+Q^{(i+1)}_l$.
The value for $Q^{(i-1)}_l$ decreases by 1 when we remove
a box $\alpha^{(i-1)}_{k+1}$, since we have
\begin{equation}
col(\alpha^{(i-1)}_{k+1})<col(\alpha^{(i)}_{k})
\leq col(\alpha^{(i)})=l
\end{equation}
(the first $<$ is by the above inequality (i),
and the next $\leq$ comes from Eq.(\ref{kantan2})).
The value for $Q^{(i)}_l$ also decreases by 1 when we remove
a box $\alpha^{(i)}_{k+1}$, since we have
\begin{equation}
col(\alpha^{(i)}_{k+1})<col(\alpha^{(i)}_{k})\leq
col(\alpha^{(i)})=l
\end{equation}
(the first $<$ is by the above inequality (ii)).
Combining these two facts, we see that the value
$Q^{(i-1)}_l-2Q^{(i)}_l$ within $p^{(i)}_l$ increases by 1 when we remove
boxes $\alpha^{(i-1)}_{k+1}$ and $\alpha^{(i)}_{k+1}$.
The value for $Q^{(i+1)}_l$ decreases by 1 or remains the same
according to whether 
\begin{equation}
col(\alpha^{(i+1)}_{k+1})\leq l
\mbox{ or }col(\alpha^{(i+1)}_{k+1})> l.
\end{equation}
However we can say that the vacancy number $p^{(i)}_l$
itself does not decrease within the region described in Eq.(\ref{kantan2}),
while we are removing boxes $\alpha^{(i-1)}_{k+1}$,
$\alpha^{(i)}_{k+1}$, and $\alpha^{(i+1)}_{k+1}$.

Next we treat the case $\alpha^{(i+1)}_{k}=\emptyset$.
Consider the region
\begin{equation}\label{kantan3}
col(\alpha^{(i)}_{k})\leq col(\alpha^{(i)}),
\end{equation}
where $\alpha^{(i)}\in\mu^{(i)}$ (see Fig.2).
We remove $\alpha^{(i-1)}_{k+1}$,
$\alpha^{(i)}_{k+1}$ and $\alpha^{(i+1)}_{k+1}$.
In this case, we can again use the above argument to get
that within the vacancy number
$p^{(i)}_l=Q^{(i-1)}_l-2Q^{(i)}_l+Q^{(i+1)}_l$,
$Q^{(i-1)}_l$ and $Q^{(i)}_l$ decrease by 1.
Thus, without any further conditions on the box $\alpha^{(i+1)}_{k+1}$,
we can deduce that the vacancy numbers $p^{(i)}_l$ do not decrease
within the region described in Eq.(\ref{kantan3}).

So far, we are discussing about the effect of the boxes
$\alpha^{(i-1)}_{k+1}$,
$\alpha^{(i)}_{k+1}$, and $\alpha^{(i+1)}_{k+1}$.
Furthermore, for some $k'>k+1$, we see that
if we remove the boxes $\alpha^{(i-1)}_{k'}$,
$\alpha^{(i)}_{k'}$, and possibly $\alpha^{(i+1)}_{k'}$,
then the vacancy numbers
for region Eq.(\ref{aida}) or region Eq.(\ref{haji}) do not decrease.
To see this, note that, by the inequalities
\begin{equation}
col(\alpha^{(j)}_{k'})<col(\alpha^{(j)}_{k+1})
\qquad (j=i-1,i),
\end{equation}
$Q^{(i-1)}_l$ and $Q^{(i)}_l$ in $p^{(i)}_l$ decrease by 1,
thus vacancy numbers do not decrease.

Combining these considerations, we see that, for each $k$,
the vacancy numbers within regions Eq.(\ref{aida}) or Eq.(\ref{haji})
do not decrease while removing the rest of the row $A$.
We summarize the results obtained in Step 2 as the following
lemma.
\begin{lemma}\label{potential}
For fixed $k$, we remove boxes
$\alpha_k^{(2)},\alpha_k^{(3)},\alpha_k^{(4)},\cdots$,
as far as possible by the KKR bijection.
Then, for each $\alpha_k^{(i)}$, we have the following
three possibilities:
\begin{enumerate}
\item[(I)] If $\alpha^{(i+1)}_k\neq\emptyset$ and also
$col(\alpha^{(i)}_{k})=col(\alpha^{(i+1)}_{k})$,
then the coquantum numbers (i.e., the vacancy numbers minus riggings
for the corresponding rows)
for the rows $\alpha^{(i)}$ of a
partition $\mu^{(i)}$ within the region
\begin{equation}
col(\alpha^{(i)}_{k})\leq col(\alpha^{(i)})
\end{equation}
do not change.
\item[(II)] If $\alpha^{(i+1)}_k\neq\emptyset$ and
$col(\alpha^{(i)}_{k})<col(\alpha^{(i+1)}_{k})$,
then the coquantum numbers for rows $\alpha^{(i)}$
within the region
\begin{equation}
col(\alpha^{(i)}_{k})\leq col(\alpha^{(i)})<
col(\alpha^{(i+1)}_{k})
\end{equation}
increase by 1 (see Fig. \ref{fig1}).
\item[(III)] If $\alpha^{(i+1)}_k=\emptyset$,
then the coquantum numbers for rows $\alpha^{(i)}$
within the region
\begin{equation}
col(\alpha^{(i)}_{k})\leq col(\alpha^{(i)})
\end{equation}
increase by 1 (see Fig. \ref{fig2}).
\end{enumerate}
For each $k$ and each partition $\mu^{(i)}$,
removal of boxes
$\alpha_k^{(2)},\alpha_k^{(3)},\alpha_k^{(4)},\cdots$,
produces a nonsingular region according to the above
(I), (II), (III), and all these regions
``accumulate'' while removing the entire row $A$
(see Fig. \ref{fig3}).
\begin{figure}
\begin{center}
\unitlength 6pt
\begin{picture}(70,12)
\put(0,0){\line(0,1){10}}
\put(0,10){\line(1,0){7}}
\qbezier(0,0)(5,0.5)(7,10)
\put(0,2){\line(1,0){3.58}}
\put(0,2.6){\line(1,0){4.1}}
\put(1,2.9){$A$}
\put(0,7){\line(1,0){6.2}}
\put(0,7.6){\line(1,0){6.4}}
\put(1.8,7.9){$B$}
\put(2,11){$\mu^{(1)}$}
\thicklines
\qbezier(5,3)(8,3)(9,4.5)
\qbezier(9,4.5)(10,6)(12,6)
\qbezier(12,6)(14,6)(19,6)
\qbezier(19,6)(23,6)(23,11)
\qbezier(4.9,2.7)(10,2.7)(29,2.7)
\qbezier(29,2.7)(35,2.7)(35,11)
\qbezier(4.9,2.3)(10,2.3)(40,2.3)
\qbezier(40,2.3)(44,2.3)(45,5.5)
\qbezier(45,5.5)(46,8.5)(48,8.5)
\qbezier(48,8.5)(52,8.5)(56,8.5)
\qbezier(56,8.5)(58,8.5)(58,11)
\thinlines
\put(12,0){\line(0,1){10}}
\put(12,10){\line(1,0){7}}
\qbezier(12,0)(17,0.5)(19,10)
\put(24,0){\line(0,1){10}}
\put(24,10){\line(1,0){7}}
\qbezier(24,0)(29,0.5)(31,10)
\put(28,6){$a_1$}
\put(36,0){\line(0,1){10}}
\put(36,10){\line(1,0){7}}
\qbezier(36,0)(41,0.5)(43,10)
\put(38,3){$a_2$}
\put(48,0){\line(0,1){10}}
\put(48,10){\line(1,0){7}}
\qbezier(48,0)(53,0.5)(55,10)
\put(60,0){\line(0,1){10}}
\put(60,10){\line(1,0){7}}
\qbezier(60,0)(65,0.5)(67,10)
\put(64.5,8.5){$a_3$}
\multiput(68,5)(0.5,0){4}{$\cdot$}
\end{picture}
\end{center}
\caption{Schematic diagram of ``accumulation''
of potential walls.
Quantum space is suppressed.
When we remove a row $A$, we obtain letters, say,
$a_1$, $a_2$, $a_3$ as image of the KKR bijection
(in the diagram, the positions
of boxes $\alpha^{(a_i-1)}_i$ are indicated
as $a_i$).
For each letter, we join all the potential walls
of Fig. \ref{fig1} and Fig. \ref{fig2} appearing at each $\mu^{(i)}$.
Then these potential walls pile up (like this diagram)
while removing the entire row $A$.}
\label{fig3}
\end{figure}
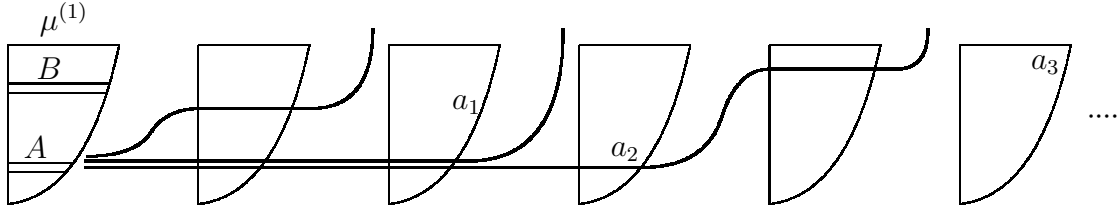
\rule{5pt}{10pt}
\end{lemma}\vspace{4mm}

\noindent
{\bf Step 3:}
We consider the consequences of Lemma \ref{potential}.
We have assumed that the rightmost part of our KKR normal ordered
product is $\cdots\otimes B\otimes A$.
We denote the width of a row $B$ as $|B|=M$.
For the sake of simplicity, we change the convention
for subscripts $k$ of $\alpha^{(2)}_k$ so that
when we remove a row $A$, then we remove the boxes
\begin{equation}
\alpha^{(2)}_{1},\,
\alpha^{(2)}_{2},\,
\cdots ,\,
\alpha^{(2)}_{m}\,
\in\mu^{(2)}|_{\leq M},
\end{equation}
in this order.
We have $col(\alpha^{(2)}_{i})>col(\alpha^{(2)}_{i+1})$,
hence $M\geq m$ holds.
We introduce one more notation.
When we remove boxes
$\alpha^{(2)}_i,\alpha^{(3)}_i,\alpha^{(4)}_i,\cdots$
as far as possible by the KKR procedure,
then we finally obtain a letter $a_i$ as an image of the KKR bijection
(i.e., $\alpha^{(a_i-1)}_i\neq\emptyset$ and
$\alpha^{(a_i)}_i=\emptyset$).
{}From the arguments in Step 2, we see that if
$\alpha^{(i)}_k\neq\emptyset$, then $\alpha^{(i)}_{k+1}\neq\emptyset$.
Interpreting this in terms of letters $a_i$, we obtain
\begin{equation}
a_1\leq a_2\leq\cdots \leq a_m.
\end{equation}

After removing a row $A$, we remove a row $B$.
Then we obtain letters $b_i$ as image of the KKR bijection,
which satisfy the inequality
\begin{equation}
b_1\leq b_2\leq\cdots \leq b_M.
\end{equation}
Then the following property holds.
\begin{proposition}\label{NYpair}
The letters $b_i$ satisfy the inequality
\begin{equation}
b_i<a_i\qquad (1\leq i\leq m).
\end{equation}
\end{proposition}
{\bf Proof.}
As a notation, when we remove the $k$th box from the right end
of a row $B$, then we remove the box $\beta^{(i)}_k$ of a partition $\mu^{(i)}$.
First, we consider the letter $b_1$.
If $\beta^{(2)}_1=\emptyset$, then $b_1=2$.
On the other hand, assuming that $m\geq 1$, i.e.,
at least there exists one box
$\alpha^{(2)}_1\in\mu^{(2)}|_{\leq M}$, then $a_1\geq 3$,
and we obtain that $b_1<a_1$ as requested.

Thus we assume that $\beta^{(2)}_1\neq\emptyset$.
We also assume that $\alpha^{(3)}_1\neq\emptyset$
(the other possibility $\alpha^{(3)}_1=\emptyset$ will
be treated later).
Then from Lemma \ref{potential} (II) we have that the rows $\beta^{(2)}$
of $\mu^{(2)}$ within the region
\begin{equation}
col(\alpha^{(2)}_1)\leq col(\beta^{(2)})<
col(\alpha^{(3)}_1)
\end{equation}
are not singular, so that $\beta^{(2)}_1$ do not fall within this region.
We have one more restriction on $\beta^{(2)}_1$.
Since $|B|=M$, we have $col(\beta^{(1)}_1)=M$,
thus $col(\beta^{(2)}_1)\geq M$ (by the definition of the KKR bijection).
On the other hand, by the definition of the present convention
(see beginning of this Step 3),
we have
$col(\alpha^{(2)}_1)\leq M.$
{}From these two inequalities we conclude
that $\beta^{(2)}_1$ must satisfy
\begin{equation}
col(\alpha^{(2)}_1)\leq col(\beta^{(2)}_1).
\end{equation}
Combining the above two restrictions on $\beta^{(2)}_1$,
we see that if $\beta^{(2)}_1\neq\emptyset$, then
\begin{equation}
col(\alpha^{(3)}_1)\leq col(\beta^{(2)}_1).
\end{equation}

Now we inductively remove boxes
$\beta^{(3)}_1$, $\beta^{(4)}_1$, $\cdots$.
Suppose that $\beta^{(i)}_1\neq\emptyset$
for some $i>2$.
First, consider the case $i<a_1-1$, i.e.,
$\alpha^{(i+1)}_1\neq\emptyset$.
Then as an induction hypothesis, we set
\begin{equation}
col(\alpha^{(i)}_1)\leq col(\beta^{(i-1)}_1).
\end{equation}
By the definition of the KKR bijection, we have
\begin{equation}
col(\beta^{(i-1)}_1)\leq col(\beta^{(i)}_1).
\end{equation}
By Lemma \ref{potential} (II)
the rows $\beta^{(i)}$ of a partition $\mu^{(i)}$
within the region
\begin{equation}
col(\alpha^{(i)}_1)\leq col(\beta^{(i)})<
col(\alpha^{(i+1)}_1)
\end{equation}
are not singular.
Therefore, if $\beta^{(i)}_1\neq\emptyset$
and $\alpha^{(i+1)}_1\neq\emptyset$, then
\begin{equation}\label{norikoe}
col(\alpha^{(i+1)}_1)\leq col(\beta^{(i)}_1).
\end{equation}
By induction, the above inequality holds
for all $i<a_1-1$.

In such a way, we remove boxes $\beta^{(i)}_1$
according to the KKR procedure.
Then, in some cases, it is possible that
$\beta^{(i)}_1=\emptyset$ for some $i\leq a_1-2$.
This means, in terms of the letters $b_1$ and $a_1$, that $b_1<a_1$.
On the contrary, it is also possible that we manage to get to
a partition $\mu^{(a_1-2)}$ and have $\beta^{(a_1-2)}_1\neq\emptyset$.
Then, as the next step, we have to consider the restrictions
imposed on the box $\beta^{(a_1-1)}_1$ (if exists).
It has to satisfy
\begin{equation}
col(\alpha^{(a_1-1)}_1)\leq
col(\beta^{(a_1-2)}_1)\leq
col(\beta^{(a_1-1)}_1),
\end{equation}
where the first $\leq$ is by Eq.(\ref{norikoe}),
and the second $\leq$ is by the definition of the KKR bijection.
We consider a next restriction. We have that
\begin{equation}
\alpha^{(a_1-1)}_1\neq\emptyset ,\quad
\alpha^{(a_1)}_1=\emptyset ,
\end{equation}
by the definition of a letter $a_1$.
Then by Lemma \ref{potential} (III) the rows
$\beta^{(a_1-1)}\in\mu^{(a_1-1)}$
within the region
\begin{equation}
col(\alpha^{(a_1-1)}_1)\leq col(\beta^{(a_1-1)}),
\end{equation} 
are not singular.
Thus we have
\begin{equation}
col(\beta^{(a_1-1)}_1)<col(\alpha^{(a_1-1)}_1),
\end{equation}
in order $\beta^{(a_1-1)}_1$ to exist.
Combining these two mutually contradicting inequalities,
we deduce that
\begin{equation}
\beta^{(a_1-1)}_1=\emptyset
\end{equation}
in any case.

To summarize, from all the above discussions we have
\begin{equation}
b_1<a_1,
\end{equation}
whenever there exist $\alpha^{(2)}_1\in\mu^{(2)}|_{\leq M}$.

Let us continue these considerations;
this time we remove boxes
$\beta^{(1)}_2$, $\beta^{(2)}_2$, $\beta^{(3)}_2$, $\cdots$.
If $\alpha^{(2)}_2\neq\emptyset$, then it has to satisfy
\begin{equation}
col(\alpha^{(2)}_2)\leq
col(\alpha^{(2)}_1)-1\leq M-1=
col(\beta^{(1)}_2).
\end{equation}

Under this setting, there is one thing that must be clarified.
\begin{lemma}
When we remove boxes
$\alpha^{(2)}_2$, $\alpha^{(3)}_2$, $\cdots$,
nonsingular regions appear on each partition
according to Lemma \ref{potential}.
Then these regions do not become singular even
after we have removed boxes
$\beta^{(1)}_1$, $\beta^{(2)}_1$, $\beta^{(3)}_1$, $\cdots$.
\end{lemma}
{\bf Proof.}
First, consider the case $\alpha^{(i+1)}_1\neq\emptyset$.
Then the rows $\alpha^{(i)}$ within the region
\begin{equation}\label{temp1}
col(\alpha^{(i)}_1)\leq
col(\alpha^{(i)})<
col(\alpha^{(i+1)}_1),
\end{equation}
of a partition $\mu^{(i)}$ are nonsingular
(by Lemma \ref{potential} (II)).
Furthermore,
since we also have $\alpha^{(i+1)}_2\neq\emptyset$ in this case,
the coquantum numbers for the rows $\alpha^{(i)}$ within the region
\begin{equation}\label{temp2}
col(\alpha^{(i)}_2)\leq
col(\alpha^{(i)})<
col(\alpha^{(i+1)}_2)
\end{equation}
of a partition $\mu^{(i)}$ (after removing the box $\alpha^{(i)}_1$)
increase by 1.
Relative locations of these two regions are described by
\begin{equation}
col(\alpha^{(i)}_2)< col(\alpha^{(i)}_1),\quad
col(\alpha^{(i+1)}_2)< col(\alpha^{(i+1)}_1).
\end{equation}
Then the following three regions of $\mu^{(i)}$ are of interest:
\begin{enumerate}
\item[(i)]   $\max \{ col(\alpha^{(i+1)}_2),col(\alpha^{(i)}_1)\}
             \leq col(\alpha^{(i)})< col(\alpha^{(i+1)}_1)$.
             The coquantum number in this region is at least 1.
\item[(ii)]  If $col(\alpha^{(i)}_1)<col(\alpha^{(i+1)}_2)$, then the region
             $col(\alpha^{(i)}_1)\leq col(\alpha^{(i)})< 
             col(\alpha^{(i+1)}_2)$ is not an empty set.
             The coquantum number of this region is at least 2.
\item[(iii)] $col(\alpha^{(i)}_2)\leq col(\alpha^{(i)})<
             \min\{ col(\alpha^{(i+1)}_2),col(\alpha^{(i)}_1)\}$.
             The coquantum number in this region is at least 1.
\end{enumerate}
Region (i) is induced by $\alpha^{(i)}_1$ and $\alpha^{(i+1)}_1$,
and region (iii) is induced by
$\alpha^{(i)}_2$ and $\alpha^{(i+1)}_2$.
Region (ii) is a superposition of these two effects.

On the other hand, from Eq.(\ref{norikoe}), we have
\begin{equation}
col(\alpha^{(i)}_1)\leq
col(\beta^{(i-1)}_1),
\end{equation}
if $\beta^{(i-1)}_1\neq\emptyset$.
This means that $Q^{(i-1)}_l$ for $l<col(\alpha^{(i)}_1)$
do not change when we remove $\beta^{(i-1)}_1$.
Similarly, when we remove boxes $\beta^{(i-1)}_1$,
$\beta^{(i)}_1$ and $\beta^{(i+1)}_1$, the vacancy number
\begin{equation}
p^{(i)}_l=Q^{(i-1)}_l-2Q^{(i)}_l+Q^{(i+1)}_l
\end{equation}
for $l<col(\alpha^{(i)}_1)$ do not change
because of the inequality
\begin{equation}
col(\alpha^{(i)}_1)\leq
col(\beta^{(i-1)}_1)\leq col(\beta^{(i)}_1)\leq
col(\beta^{(i+1)}_1).
\end{equation}
As a result, region (iii) in the above do not become singular
after removing $\beta^{(i-1)}_1$, $\beta^{(i)}_1$,
and $\beta^{(i+1)}_1$.

The coquantum number for region (ii) above might decrease by 1
when we remove $\beta^{(i-1)}_1$,
$\beta^{(i)}_1$, and $\beta^{(i+1)}_1$.
For example, if
\begin{equation}
col(\beta^{(i-1)}_1)\leq
col(\alpha^{(i+1)}_2)<col(\beta^{(i)}_1),
\end{equation}
then $Q^{(i-1)}_l$ decrease by 1, and $Q^{(i)}_l$ and
$Q^{(i+1)}_l$ do not change
(where $l=col(\alpha^{(i)})$ for a box $\alpha^{(i)}$
within the region (ii) above).
However, the coquantum numbers for region (ii) are more than 2,
thus region (ii) also does not become singular.
After all, we see that nonsingular region Eq.(\ref{temp2})
(=(ii) $\cup$ (iii) in the above classification)
does not become singular even if we remove boxes
$\beta^{(i-1)}_1$, $\beta^{(i)}_1$, and $\beta^{(i+1)}_1$.

The case $\alpha^{(i+1)}_1=\emptyset$ is almost similar.
We use Lemma \ref{potential} (III) and
$col(\alpha^{(i)}_1)\leq col(\beta^{(i-1)}_1)$ (from Eq.(\ref{norikoe})).
In this case, the following two regions of $\mu^{(i)}$ are of interest:
\begin{enumerate}
\item[(i)$'$] $col(\alpha^{(i)}_2)\leq col(\alpha^{(i)})
              <col(\alpha^{(i)}_1)$ if $\alpha^{(i+1)}_2=\emptyset$, or
              $col(\alpha^{(i)}_2)\leq col(\alpha^{(i)})
              <\min\{ col(\alpha^{(i+1)}_2),col(\alpha^{(i)}_1)\}$
              if $\alpha^{(i+1)}_2\neq\emptyset$. The coquantum numbers in these
              regions are at least 1.
\item[(ii)$'$]$col(\alpha^{(i)}_1)\leq col(\alpha^{(i)})$ if
              $\alpha^{(i+1)}_2=\emptyset$, or 
              $col(\alpha^{(i)}_1)\leq col(\alpha^{(i)})
              <col(\alpha^{(i+1)}_2)$ if $\alpha^{(i+1)}_2\neq\emptyset$ and
              $col(\alpha^{(i)}_1)<col(\alpha^{(i+1)}_2)$.
              The coquantum numbers for these regions are at least 2.
\end{enumerate}
Region (i)$'$ is induced by $\alpha^{(i)}_2$, and
region (ii)$'$ is induced by both $\alpha^{(i)}_1$ and $\alpha^{(i)}_2$.
Since $col(\alpha^{(i)}_1)\leq col(\beta^{(i-1)}_1)$,
region (i)$'$ does not become singular, and since coquantum numbers of
region (ii)$'$ are at least 2, it also does not become singular.

This completes the proof of the lemma.
\rule{5pt}{10pt}\vspace{4mm}

Keeping this lemma in mind, let us return to the proof of the proposition.
We remove boxes
$\beta^{(2)}_2$, $\beta^{(3)}_2$, $\beta^{(4)}_2$, $\cdots$
as far as possible.
When we removed boxes $\alpha^{(2)}_2$, $\alpha^{(3)}_2$, $\cdots$
there are regions of partitions whose coquantum numbers
increased according to Lemma \ref{potential}.
Before the above lemma, we have shown that
$col(\alpha^{(2)}_2)\leq col(\beta^{(1)}_2)$.
Thus we can apply the argument which was used when we removed boxes
$\beta^{(1)}_1$, $\beta^{(2)}_1$, $\beta^{(3)}_1$, $\cdots$
to get
\begin{equation}
b_2<a_2.
\end{equation}

We can apply the same argument to the remaining letters
$a_3$, $a_4$, $\cdots$ and get
\begin{equation}
b_i<a_i,\qquad
(1\leq i\leq m).
\end{equation}

This completes the proof of the proposition.
\rule{5pt}{10pt}\vspace{4mm}

\noindent
{\bf Step 4:}
In this and the following steps, we calculate the unwinding number of $B\otimes A$
based on the above considerations.
First of all, we make the following distinctions.
The row $\alpha$ of $\mu^{(2)}$, which is removed when we remove a row $A$,
is the shortest row among the rows of $\mu^{(2)}$
whose widths $w$ satisfy $w\geq M$
before we remove a row $A$.
When we remove a row $A$, the row $\alpha$ is removed to be the row $\alpha '$.
Then there are the following three possibilities:
\begin{enumerate}
\item[(a)] $col(\alpha ')>M$,
\item[(b)] $col(\alpha ')\leq M$,
\item[(c)] there is no such a row $\alpha$, i.e., all boxes of $\mu^{(2)}$
           which are removed with a row $A$ are elements of $\mu^{(2)}|_{\leq M}$.
\end{enumerate}
In this step, we treat case (a).

We continue to use the notation of Step 3;
when we remove a row $A$, then the elements of $\mu^{(2)}|_{\leq M}$
\begin{equation}
\alpha^{(2)}_1, \alpha^{(2)}_2, \cdots ,\alpha^{(2)}_m
\end{equation}
are removed in this order.
We consider the box
\begin{equation}
\alpha^{(2)}_0\in\mu^{(2)}|_{>M},
\end{equation}
which is the last box among all the boxes of the row $\alpha '$ 
that are removed with a row $A$.
By the KKR procedure, we remove boxes $\alpha^{(2)}_0$,
$\alpha^{(3)}_0$, $\alpha^{(4)}_0$, $\cdots$
as far as possible and eventually obtain a letter $a_0$.
In other words, we have $\alpha^{(a_0-1)}_0\neq\emptyset$
and $\alpha^{(a_0)}_0=\emptyset$.

After removing boxes $\alpha^{(2)}_0$,
$\alpha^{(3)}_0$, $\cdots$, the remaining rows $\alpha^{(2)}_0-1$,
$\alpha^{(3)}_0-1$, $\cdots$ are made to be singular.
Then the simplest case is as follows.
We assume that these singular rows remain to be singular
even after a row $A$ is entirely removed.

We remove the rightmost box of a row $B$, i.e.,
box $\beta^{(1)}_1\in\mu^{(1)}$.
Then it satisfies $col(\beta^{(1)}_1)=M$.
In the next partition $\mu^{(2)}$,
the row $\alpha^{(2)}_0-1$ is singular, and its width is
$col(\alpha^{(2)}_0-1)\geq M$.
Thus, in one case, we can remove the boxes $\beta^{(1)}_1$,
$\alpha^{(2)}_0-1$, $\alpha^{(3)}_0-1$, $\cdots$,
$\alpha^{(a_0-1)}_0-1$, $\cdots$,
or, in the other case, we remove the boxes
$\beta^{(1)}_1$, ${\alpha '}^{(2)}_0$, ${\alpha '}^{(3)}_0$, $\cdots$,
which satisfy $col({\alpha '}^{(i)}_0)< col(\alpha^{(i)}_0-1)$.
In the former case, we have
\begin{equation}
b_1\geq a_0.
\end{equation}
In the latter case, if $i\leq a_0-2$, then we always have singular rows
$\alpha^{(i+1)}_0-1$ which satisfy 
$col({\alpha '}^{(i)}_0)<col(\alpha^{(i+1)}_0-1)$.
Thus we can remove the boxes $\alpha^{(i+1)}_0-1$ if necessary,
so that we have
\begin{equation}
b_1\geq a_0.
\end{equation}
Hence we obtain
\begin{equation}
b_i\geq a_0\qquad
(1\leq i\leq M),
\end{equation}
because of the inequalities $b_{i+1}\geq b_i$.

The next simplest case is as follows.
After removing a row $A$, the rows
\begin{equation}
\alpha^{(2)}_0-1, \alpha^{(3)}_0-1, \cdots ,\alpha^{(i-1)}_0-1,
\end{equation}
remain singular; on the other hand,
the coquantum number of the row $\alpha^{(i)}_0-1$ for some
$i\leq a_0$ becomes 1.

Since the coquantum number of the row $\alpha^{(i)}_0-1$ is increased by 1,
we can deduce the following two possibilities
by use of Lemma \ref{potential};
the box $\alpha^{(i)}_0-1$ is within either
\begin{equation}\label{step4_1}
col(\alpha^{(i)}_1)\leq col(\alpha^{(i)}_0-1)
<col(\alpha^{(i+1)}_1)
\end{equation}
if $\alpha^{(i+1)}_1\neq\emptyset$, or
\begin{equation}
col(\alpha^{(i)}_1)\leq col(\alpha^{(i)}_0-1)
\end{equation}
if $\alpha^{(i+1)}_1=\emptyset$.
(In view of $col(\alpha^{(i)}_1)<col(\alpha^{(i)}_0)$,
we have yet another possibility
\begin{equation}
col(\alpha^{(i+1)}_1)\leq col(\alpha^{(i)}_0-1)
\end{equation}
when $\alpha^{(i+1)}_1\neq\emptyset$.
However we need not take it into consideration, since,
in such a situation, the vacancy number of the row
$\alpha^{(i)}_0-1$ does not change.)

First, assume that $\alpha^{(i+1)}_1\neq\emptyset$.
Under these settings, we further assume that the rows
$\alpha^{(i+1)}_1-1$, $\alpha^{(i+2)}_1-1$, $\cdots$,
$\alpha^{(a_1-1)}_1-1$ remain singular even
after removing a row $A$.
Then we shall show that the inequality
\begin{equation}
b_2\geq a_0
\end{equation}
holds in this case.
We shall generalize these arguments later.

When we begin to remove a row $B$, then we remove
$\beta^{(1)}_1$, $\beta^{(2)}_1$, $\beta^{(3)}_1$, $\cdots$
as far as possible and obtain a letter $b_1$
as the image of the KKR bijection.
By the above assumption, the rows
$\alpha^{(2)}_0-1$, $\alpha^{(3)}_0-1$, $\cdots$, $\alpha^{(i-1)}_0-1$
are singular, and $col(\beta^{(1)}_1)\leq col(\alpha^{(2)}_0-1)$.
Therefore, at least we have 
\begin{equation}
\beta^{(i-1)}_1\neq\emptyset ,
\end{equation}
i.e.,
$b_1\geq i,$
and we also have
\begin{equation}
col(\beta^{(j)}_1)\leq col(\alpha^{(j)}_0)-1\qquad
(1\leq j\leq i-1).
\end{equation}
We consider in turn the possible states of $\beta^{(i)}_1$.
When $\beta^{(i)}_1=\emptyset$, then the coquantum number of
the row $\alpha^{(i)}_0-1$ decreases by 1 so that it becomes singular.
Next, when $col(\beta^{(i)}_1)>col(\alpha^{(i)}_0)-1$,
then the coquantum number of
the row $\alpha^{(i)}_0-1$ decreases by 1 independently of the position
of $\beta^{(i+1)}_1$.
In these two cases, if we remove $\beta^{(1)}_2$, at least
we can remove
\begin{equation}
\alpha^{(2)}_0-1, \cdots ,\alpha^{(i)}_0-1,
\alpha^{(i+1)}_1-1, \cdots ,\alpha^{(a_1-1)}_1-1,
\end{equation}
or, in terms of letters $a_i$ and $b_i$, we have
\begin{equation}
b_2\geq a_1\geq a_0.
\end{equation}
On the other hand, consider the case
$col(\beta^{(i)}_1)\leq col(\alpha^{(i)}_0)-1$.
Then from the above discussion (see Eq.(\ref{step4_1})),
we already have the restriction
\begin{equation}
col(\alpha^{(i)}_0)-1\leq col(\alpha^{(i+1)}_1-1),
\end{equation}
thus we can remove $\alpha^{(i+1)}_1-1$ as $\beta^{(i+1)}_1$.
Therefore we deduce that $b_1\geq a_1$, i.e.,
\begin{equation}
b_2\geq b_1\geq a_1\geq a_0.
\end{equation}

The case $\alpha^{(i+1)}_1=\emptyset$ is similar.
In this case, we also have $\beta^{(i-1)}_1\neq\emptyset$.
Since, in this case, $i=a_1-1$ if $\beta^{(i)}_1\neq\emptyset$, then
$b_1\geq a_1$, i.e.,
\begin{equation}
b_2\geq b_1\geq a_1\geq a_0.
\end{equation}
On the other hand, if $\beta^{(i)}_1=\emptyset$,
then from the inequality
\begin{equation}
col(\beta^{(i-1)}_1)\leq
col(\alpha^{(i-1)}_0)-1\leq
col(\alpha^{(i)}_0)-1
\end{equation}
we have that the coquantum number of the
row $\alpha^{(i)}_0-1$ decreases by 1
and it becomes singular.
Thus we conclude that
\begin{equation}
b_2\geq a_0.
\end{equation}

In the above discussion, we have shown that $b_2\geq a_0$
under some restriction. We can generalize the arguments
as follows.
\begin{enumerate}
\item[(i)]  In the above arguments, we have assumed that the rows
$\alpha^{(i)}_1-1$, $\alpha^{(i+1)}_1-1$, $\cdots$ remain
singular even if we remove a row $A$.
To generalize it, we consider as follows.
The rows $\alpha^{(1)}_0-1$, $\cdots$, $\alpha^{(i_1-1)}_0-1$
are singular, however, the coquantum number of the row $\alpha^{(i_1)}_0-1$
becomes 1 as in the above arguments.
Next, $\alpha^{(i_1+1)}_1-1$, $\cdots$, $\alpha^{(i_2-1)}_1-1$
are singular, however, the coquantum number of $\alpha^{(i_2)}_1-1$
becomes 1 because of removal of $\alpha^{(i_2)}_2$ and others,
$\cdots$, $\alpha^{(i_{k-1}+1)}_{k-1}-1$,
$\cdots$, $\alpha^{(i_k-1)}_{k-1}-1$
are singular, however, the coquantum number of $\alpha^{(i_k)}_{k-1}-1$
becomes 1 because of removal of $\alpha^{(i_k)}_k$ and others,
and rows $\alpha^{(i_k+1)}_k-1$, $\cdots$, $\alpha^{(a_k-1)}_k-1$
remain to be singular.
Then, by applying the above arguments to each step,
we see that if we remove at least $k$ boxes from a row $B$,
then the sequence $\beta^{(1)}_{k+1}$, $\beta^{(2)}_{k+1}$, $\cdots$
satisfies $\beta^{(i_k+1)}_{k+1}\neq\emptyset$; therefore we obtain
\begin{equation}
b_{k+1}\geq a_k\geq a_0.
\end{equation}
\item[(ii)] On the other hand, it is possible that, after removing
a row $A$, the rows $\alpha^{(1)}_0-1$, $\alpha^{(2)}_0-1$, $\cdots$,
$\alpha^{(i-1)}_0-1$ remain singular, but the coquantum number of
the row $\alpha^{(i)}_0-1$ is $k$.
In this case, we can also apply the above arguments to show that
the coquantum number of the row $\alpha^{(i)}_0-1$ is raised by $k$ 
because of
\begin{equation}
\alpha^{(i)}_1, \alpha^{(i)}_2, \cdots ,\alpha^{(i)}_k
\end{equation}
and the adjacent ones.
In this case, if we remove at least $k$ boxes from a row $B$,
then the row $\alpha^{(i)}_0-1$ becomes singular.
Thus, if we remove $\beta^{(1)}_{k+1}$, then we can remove
$\alpha^{(i+1)}_k-1$, $\alpha^{(i+2)}_k-1$, $\cdots$.
The fundamental case is that the rows $\alpha^{(j)}_k-1$
$(j\geq i+1)$ remain singular, and, in this case, we have
\begin{equation}
b_{k+1}\geq a_k\geq a_0.
\end{equation}
\end{enumerate}
We can combine the above (i) and (ii) to treat the general case.
Especially, we notice that the relevant boxes are
\begin{equation}
\alpha^{(2)}_1,\cdots ,\alpha^{(2)}_m\in\mu^{(2)}|_{\leq M},
\end{equation}
and at least we have
\begin{equation}
b_i\geq a_0\qquad (i\geq m+1).
\end{equation}

Summarizing the above arguments of Step 4, we obtain the following:
\begin{lemma}\label{oshidashi}
In the above setting, we have that
\begin{equation}
b_i\geq a_0\qquad (i\geq m+1),
\end{equation}
where $m$ is the number of boxes removed from
$\mu^{(2)}|_{\leq M}$ when we remove the row $A$.
\rule{5pt}{10pt}
\end{lemma}

Using Lemma \ref{oshidashi},
one can derive the unwinding number of the tensor product $B\otimes A$.
By Proposition \ref{NYpair} we can connect each $b_i$ and $a_i$
$(1\leq i\leq m)$ as an unwinding pair.
On the other hand, we have
\begin{equation}
a_i\geq a_0\qquad (1\leq i\leq m),
\end{equation}
that is, there are only $m$ letters $a_i$ greater than $a_0$,
and we know that all these letters are already connected
with $b_1$, $\cdots$, $b_m$.
By Lemma \ref{oshidashi}, $b_{m+1}$, $\cdots$, $b_M$ are greater
than $a_0$; so they cannot be connected with the rest of the letters
produced by the row $A$.

As a result, if the number of letters removed from
$\mu^{(2)}|_{\leq M}$ while removing the row $A$ is $m$
and if condition (a) at the beginning of Step 4 is fulfilled, then
\begin{equation}
\mbox{the unwinding number of }B\otimes A=m,
\end{equation}
as desired.\vspace{4mm}

\noindent
{\bf Step 5:}
We considered case (b) at the beginning of Step 4.
In this case, we can apply almost similar arguments of Step 4.
Suppose that the row $\alpha '$ which appeared in case (b) satisfies
\begin{equation}
M-col(\alpha ')=l.
\end{equation}
Then the number of $\alpha^{(2)}_i$ within $\mu^{(2)}|_{<col(\alpha ')}$
is $m-l$.

In removing the row $B$, if we remove $l$ boxes from the row $B$,
then Lemma \ref{oshidashi} becomes applicable.
As a notation,
if we remove a box $\alpha ' +1$
(right adjacent of the $\alpha '$) while removing
the row $A$, then we obtain a letter $a_0'$.
By Lemma \ref{oshidashi} we have
\begin{equation}
b_i\geq a_0' \qquad (i\geq m-l+1).
\end{equation}
On the other hand, from Proposition \ref{NYpair} we have
\begin{equation}
b_i<a_i,\qquad (1\leq i\leq m).
\end{equation}
By the definition of $\alpha^{(2)}_i\in\mu^{(2)}|_{\leq M}$,
we have
\begin{equation}
a_0'\geq a_1.
\end{equation}
Then by an argument similar to that at the end of Step 4,
we have
\begin{equation}
\mbox{the unwinding number of }B\otimes A=m
\end{equation}
for case (b).\vspace{4mm}

\noindent
{\bf Step 6:}
In this step, we treat condition (c) at the beginning of Step 4.
When we remove the row $A$, we remove
\begin{equation}
\alpha^{(2)}_1,\cdots ,\alpha^{(2)}_m,
\end{equation}
and, in this case, all these boxes are elements of $\mu^{(2)}|_{\leq M}$.
If $m=0$, then $A=\fbox{$2^L$}$, where we set $|A|=L$,
so that the unwinding number of $B\otimes A$ is always equal to 0,
as was to be shown.

We assume that $m\neq 0$.
We denote the number of letters 2 in tableau $A$ as
\begin{equation}
L-m=:t.
\end{equation}
These $t$ letters 2 do not contribute to the unwinding number
of $B\otimes A$.
{}From Proposition \ref{NYpair} we have
\begin{equation}
b_i<a_i\qquad (1\leq i\leq m).
\end{equation}
Since $t+m=L$, we have already checked all letters in $A$.
Thus we also have
\begin{equation}
\mbox{the unwinding number of }B\otimes A=m
\end{equation}
in this case (c).\vspace{4mm}

Now we have shown that cases (a), (b), and (c) appearing in Step 4
all satisfy Theorem \ref{unwinding}.
Hence the proof of Theorem is finished.
\rule{5pt}{10pt}

\vspace{6mm}

\noindent
{\bf Acknowledgements:}
The author would like to thank Professor M. Wadati
for warm encouragements during the study.
He is also grateful to Professors A. Kuniba,
S. Naito, M. Okado, D. Sagaki,
A. Schilling, M. Shimozono,
T. Takagi, and Y. Yamada for discussions
and to referees for numerous suggestions about this manuscript.
He is a research fellow of the 
Japan Society for the Promotion of Science.

\end{document}